\documentclass[a4paper,11pt, oneside]{amsart}
\usepackage[left=3cm,right=3cm]{geometry}

\usepackage[T1]{fontenc}
\usepackage[english]{babel}
\usepackage{enumitem}
\usepackage{young}
\usepackage{amsthm, amssymb, amsfonts, amsmath, mathrsfs}
\usepackage[all]{xy}
\usepackage{ytableau}
\usepackage{caption}
\usepackage{url}
\usepackage{tikz}
\usepackage{hyperref}
\usepackage{comment}

\usepackage{xcolor}

\usepackage{array}
\newcommand{\PreserveBackslash}[1]{\let\temp=\\#1\let\\=\temp}
\newcolumntype{C}[1]{>{\PreserveBackslash\centering}p{#1}}
\newcolumntype{R}[1]{>{\PreserveBackslash\raggedleft}p{#1}}
\newcolumntype{L}[1]{>{\PreserveBackslash\raggedright}p{#1}}

\newcounter{stepcounter}
\theoremstyle{plain}

\newtheorem{thm}{Theorem}[section]

\newtheorem{lem}[thm]{Lemma}
\newtheorem{prop}[thm]{Proposition}
\newtheorem{cor}[thm]{Corollary}

\theoremstyle{definition}

\newtheorem{eg}[thm]{Example}
\newtheorem{defn}[thm]{Definition}

\newtheorem{remark}[thm]{Remark}

\date{}

\newcommand\bit{\begin{itemize}}
\newcommand\eit{\end{itemize}}
\newcommand\bet{\begin{enumerate}}
\newcommand\eet{\end{enumerate}}
\newcommand\ed{\end{document}}

\DeclareFontFamily{U}{mathx}{\hyphenchar\font45}
\DeclareFontShape{U}{mathx}{m}{n}{
      <5> <6> <7> <8> <9> <10>
      <10.95> <12> <14.4> <17.28> <20.74> <24.88>
      mathx10
      }{}
\DeclareSymbolFont{mathx}{U}{mathx}{m}{n}
\DeclareFontSubstitution{U}{mathx}{m}{n}
\DeclareMathAccent{\widecheck}{0}{mathx}{"71}
\DeclareMathAccent{\wideparen}{0}{mathx}{"75}

\newcommand{\e}{\varepsilon}
\newcommand\Om{\Omega}
\newcommand\del{\partial}
\newcommand\adel{\ol{\partial}}
\newcommand\DEL{\Delta}

\newcommand\bC{{\mathbb C}}

\newcommand\bZ{{\mathbb Z}}

\newcommand\A{{\mathcal{A}}}

\newcommand\EE{{\mathcal E}}
\newcommand\F{{\mathcal F}}

\newcommand{\OO}{\mathcal{O}}

\newcommand\Ad{\mathrm{Ad}}

\newcommand\can{\mathrm{can}}

\newcommand\co{\mathrm{co}}

\newcommand\exd{\mathrm{d}}

\newcommand\haar{\mathrm{\bf h}}

\newcommand\unit{\mathrm{U}}
\newcommand\counit{\mathrm{C}}

\newcommand\id{\mathrm{id}}
\newcommand\proj{\mathrm{proj}}

\newcommand\spn{\mathrm{span}}
\newcommand\sstab{\mathrm{SSTab}}

\newcommand\hol{^{(1,0)}}
\newcommand\ahol{^{(0,1)}}

\newcommand\inv{^{-1}}

\newcommand\by{\times}

\newcommand\oby{\otimes}

\newcommand\wed{\wedge}
\newcommand\sseq{\subseteq}
\newcommand\tl{\triangleleft}

\def\qbinom#1#2{\ensuremath{\left[\kern-.3em\left[\genfrac{}{}{0pt}{}{#1}{#2}\right]\kern-.3em\right]_q}}

\newcommand\ol{\overline}

\newcommand\bs{\backslash}
\newcommand\mto{\mapsto}

\DeclareMathOperator{\dt}{det}

\usepackage{tikz}
\usepackage{tikz-cd}
\usetikzlibrary{decorations.pathreplacing}

\usetikzlibrary{arrows}
\usepackage[english]{babel}

\def\cprime{\/{\mathsurround=0pt$'\!\!$}}

\usepackage{mathtools}

\setenumerate[1]{label={\arabic*.}}
\title{A Borel--Weil Theorem for the Quantum Grassmannians}

\author[R. \'O Buachalla]{R\'eamonn \'O Buachalla}
\address{Mathematical Institute of Charles University, Sokolovsk\'a 83, Prague, Czech Republic} \email{obuachalla@karlin.mff.cuni.cz}

\author[A. Carotenuto]{Alessandro Carotenuto}
\address{Institute of Mathematics, Czech Academy of Sciences, \v{Z}itn\'a 25, 115 67 Prague, Czech Republic}
\email{carotenuto@math.cas.cz, acaroten91@gmail.com}

\author[C. Mrozinski]{Colin Mrozinski}
\address{Instytut Matematyczny, Polskiej Akademii Nauk, ul. \'Sniadeckich 8, 00-656 Warszawa, Poland}
\email{colin.mrozinski@gmail.com}

\thanks{}


\thanks{AC is supported by the GA\v{C}R project 20-17488Y and \mbox{RVO: 67985840}. R\'OB was supported by the international cofunded project 3038/7.PR/2014/2 in the years 2014-2018, and by the grant FP7-PEOPLE-2012-COFUND-600415, he acknowledges FNRS support through  a postdoctoral fellowship within the framework of the MIS Grant ``Antipode'' grant number F.4502.18, he acknowledges GA\v{C}R support through the framework of the grant GA19-06357S, and finally support of the Charles University PRIMUS grant PRIMUS/21/SCI/026.}

\begin{document}

\maketitle

\begin{abstract}
We establish a noncommutative generalisation of  the  Borel--Weil theorem for the Heckenberger--Kolb calculi of the quantum Grassmannians. The result is formulated in the framework of quantum principal bundles and noncommutative complex structures, and generalises previous work of a number of authors on quantum projective space. As a direct consequence we get a novel noncommutative differential geometric presentation of the twisted Grassmannian coordinate ring  studied in noncommutative projective geometry. A number of applications to the noncommutative K\"ahler geometry of the quantum Grassmannians are also given. \end{abstract}

\tableofcontents

\section{Introduction}


The Borel--Weil theorem \cite{serre} is  an elegant  geometric  procedure for constructing all unitary irreducible representations of a compact Lie group. The construction realises each representation as the space of holomorphic sections of a  line bundle over a flag manifold. It is a highly influential result in the representation theory  of Lie groups and since the discovery of quantum groups has inspired a number of noncommutative generalisations. Important examples include the derived functor approach of Andersen, Polo, Wen \cite{APW}, 
the quantum coordinate algebra approaches of Parshall, Wang \cite{PW}, and  Mimachi, Noumi, Yamada \cite{ noumi0,NoumiUn}, the approach of Biedenharn, Lohe \cite{BL} based on q-bosons, the coherent state approach of Jur\v{c}o, \v{S}{\v{t}}ov\'i\v{c}ek, \cite{JuSt}, 
 and the compact quantum group approach of Gover, Zhang \cite{GZ}. Moreover, in \cite{KremnitzerBack,KremnitzerBakROU} Kremnitzer and Backelin generalised Beilinson--Bernstein localisation (itself a far-reaching generalisation of the Borel--Weil theorem) to the setting of quantum flag manifolds.


In the classical setting the holomorphic sections of a  line bundle over a flag manifold are the same as its Borel invariant sections. Very roughly speaking, the above works generalise the Borel invariant description without introducing any formal noncommutative notion of  holomorphicity. 
In recent years, however, the study of differential graded algebras over the quantum flag manifolds $\OO_q(G/L_S)$ has yielded a much better understanding of their noncommutative complex geometry. Subsequent work on the Borel--Weil theorem for quantum groups has used the notion of a complex structure on a differential $*$-calculus  to generalise the Koszul--Malgrange presentation of holomorphic vector bundles  \cite{KM}. This direction of research was initiated by Majid in his influential paper  on  the Podle\'s sphere \cite{Maj}. It was continued by Khalkhali, Landi, van Suijlekom, and Moatadelro in \cite{KLVSCP1,KKCP2, KKCPN} where the definitions of  complex structure and noncommutative holomorphic vector bundle were introduced and the family of examples extended to include quantum projective space $\OO_q(\mathbb{CP}^{n})$. The same notion of holomorphic structure would later appear independently in the work of Beggs and Smith on noncommutative coherent sheaves \cite{BS}.


The differential calculi used in the above works are those identified by Heckenberger and Kolb in their remarkable classification of calculi over  the irreducible quantum flag manifolds \cite{HK}. This is one of the most important results in the study of the noncommutative geometry of quantum groups, and a touchstone for future theory. In the last years remarkable progress has been made in our understanding of the Heckenberger--Kolb calculi, in particular in our understanding of their noncommutative complex and K\"ahler geometry, see for example \cite{MarcoConj,MMF2,MMF3}. Of particular relevance to this paper is the recent proof that every homogeneous vector module over an irreducible quantum flag manifold $\OO_q(G/L_m)$ admits a unique left $\OO_q(G)$-covariant holomorphic structure \cite{HVBQFM}. For the special case of quantum projective space $\OO_q(SU_n/U_{n-1})$ this reduces to the holomorphic structures explicitly constructed in \cite{KLVSCP1,KKCP2,KKCPN,Maj}, giving them a conceptual presentation. Thus it is natural to ask if the $q$-deformed Borel--Weil theorem for quantum projective space extends to all the irreducible quantum flag manifolds. In this paper we show that this is indeed the case for the $A$-series examples, which is to say the quantum Grassmannians.


We approach the problem using the framework of quantum principal bundles. This means that the paper is closer in form to \cite{KLVSCP1,Maj} rather than \cite{KKCP2, KKCPN}, where D\c{a}browski and D'Andrea's spectral triple presentation of the quantum projective space calculus was used.
We start by taking the direct sum of all the line bundle over the quantum Grassmannian $\OO_q(\mathrm{Gr}_{n,m})$. This forms a $\OO(U_1)$-principal comodule algebra which we call the quantum Grassmann sphere. For quantum projective space it reduces to the odd-dimensional quantum sphere, and for the special case of the quantum projective line, we recover the celebrated quantum Hopf fibration over the Podle\'s sphere.
In order to construct a differential structure for the Hopf--Galois extension $\OO_q(\mathrm{Gr}_{n,m}) \hookrightarrow \OO_q(S^{n,m})$, we generalise the approach introduced in \cite{MMF1}. Explicitly, we exploit the well-known connection between coquasi-triangular structures and bicovariant differential first-order calculi \cite{BeggsMajid:Leabh,JurcoBC,KSLeabh} to construct a quantum principal bundle structure for the Hopf--Galois extension 
$\OO_q(\mathrm{Gr}_{n,m}) \hookrightarrow \OO_q(S^{n,m})$. This allows us to induce the holomorphic structures of the line bundles  $\EE_k$ from a principal connection for the bundle, and it is this picture that allows us to prove our quantum generalisation of the Borel--Weil theorem.


Throughout the paper we adopt a quantum coordinate algebra approach, using the FRT-construction of $\OO_q(SU_n)$ and its coquasitriangular structure. Indeed, with a view to minimising the necessary preliminaries, the quantised enveloping algebra $U_q(\mathfrak{sl}_n)$ is discussed only in an expository appendix. As such, for the comodule classification and the Peter--Weyl decomposition of $\OO_q(SU_n)$ we follow the quantum minor formulation of Mimachi, Noumi, and Yamada. This approach fits particularly well with the quantum Grassmannians: the algebra $\OO_q(S^{n,m})$ is generated as a $*$-algebra by certain quantum minors of length $n-m$, and endowing the minors and their conjugates with degree $1$ and $-1$ respectively recovers the $\mathbb{Z}$-grading of $\OO_q(S^{n,m})$. Moreover, the holomorphic elements of each line bundle $\EE_k$ are precisely the degree-$k$ elements of the homogeneous coordinate ring $S_q(\mathrm{Gr}_{n,m})$, giving a geometric realisation of the Borel coinvariant presentation of the corresponding irreducible $\OO_q(SU_n)$-comodule.


An important motivation of the paper is to further explore the connections between the quantum coordinate algebras $\OO_q(G/L_S)$ and their twisted homogeneous coordinate ring counterparts $S_q(G/L_S)$ in noncommutative projective geometry \cite{Avdb,Soib}. For the special case of the quantum Grassmannians, these rings are  deformations of the Pl\"ucker embedding homogeneous coordinate ring and are important examples in the theory of quantum cluster algebras \cite{KKCP2,KKCPN}. 
The quantum Grassmannian Borel--Weil theorem gives us a direct $q$-deformation of the classical ample line bundle presentation of  $S_q(\mathrm{Gr}_{n,m})$. This directly generalising the work of \cite{KLVSCP1,KKCP2, KKCPN} for quantum projective space,
and gives us an important point of contact between  noncommutative differential geometry and noncommutative projective geometry (see \cite{BS} for a discussion of what a noncommutative generalisation of the classical GAGA correspondence might look like). Another important application of the quantum Grassmannian Borel--Weil theorem is to the study of the noncommutative K\"ahler geometry of $\OO_q(\mathrm{Gr}_{n,m})$. The given cohomological information allows us to identify which line bundles over $\OO_q(\mathrm{Gr}_{n,m})$ are positive and which are negative. Moreover, it allows us to conclude that twisting the Dolbeault--Dirac operator of  $\OO_q(\mathrm{Gr}_{n,m})$ by a negative line bundle produces a Fredholm operator, which is to say, it allows us to conclude analytic behaviour from purely geometric information.


The paper naturally leads to a number of future projects. In particular, we highlight the task of extending the Borel--Weil theorem to all irreducible quantum flag manifolds, as well as the more general question of calculating the cohomology of general covariant holomorphic modules \cite[Remark 4.7]{HVBQFM}.


The paper is organised as follows:  In Section 2, we recall a number of basic Hopf algebraic constructions: Hopf--Galois extensions, principal comodule algebras, quantum homogeneous spaces, and Takeuchi's equivalence for relative Hopf modules. We recall how such structures interact with differential calculi, focusing on the theory of quantum principal bundles.  We also  present some more recent material  about  complex structures, and holomorphic structures.  

In section 3 we introduce the notion of a principal pair, a formal structure that abstracts the algebraic properties of the quantum Grassmann sphere $\OO_q(\mathrm{S}^{n,m})$.  This is a special type of Hopf--Galois extension $B = P^{\co(T)}$, where $P$ and $B$ are assumed to be quantum  homogeneous $A$-spaces, for some Hopf algebra $A$. We then introduce a robust set of results for constructing quantum principal bundles and principal connections from principal pairs. 

In Section 4, we recall the FRT-construction of $\OO_q(SU_n)$, the construction of the quantum Grassmannians $\OO_q(\mathrm{Gr}_{n,m})$, and the quantum Grassmann spheres $\OO_q(S^{n,m})$, as quantum homogeneous spaces, and recall a set of generators for both spaces. Moreover, we present the $\OO(U_1)$-Hopf--Galois extension $\OO_q(\mathrm{Gr}_{n,m}) \hookrightarrow \OO_q(S^{n,m})$ as a principal pair, and construct an explicit strong principal connection  for the bundle.

In Sections 4 and 5, we use the coquasitriangular quantum Killing form, and its associated bicovariant calculus, to construct a calculus  $\Omega^1_q(SU_n,m)$ which restricts to the Heckenberger--Kolb calculus on $\OO_q(\mathrm{Gr}_{n,m})$. Moreover, using the general results of \textsection 3, we show that $\Omega^1_q(SU_n,m)$ induces a quantum principal bundle structure on the Hopf--Galois extensions $\OO_q(\mathrm{Gr}_{n,m}) \hookrightarrow \OO_q(S^{n,m})$. The universal principal connection introduced in Section 3 is then shown to restrict to a principal connection for the bundle.

In Section 6, we prove the main result of the paper, which describes the holomorphic sections of the holomorphic line bundles $\EE_k$ introduced in \cite{DOKSS}. 

\begin{thm}[Borel--Weil]
Defining the space of holomorphic sections of the line bundle $\EE_k$ to be the $\OO_q(SU_n)$-module
\begin{align*}
H^0_{\adel}(\EE_{k}) := \ker(\adel_{\EE_k}:\EE_k \to \Omega^{(0,1)} \otimes_{\OO_q(\mathrm{Gr}_{n,m})} \EE_k),
\end{align*}
it holds that
\begin{align*}
H^0_{\adel}(\EE_{k}) =  V_{k\varpi_{n-m}}, & & H^0_{\adel}(\EE_{-k}) = 0, 
\end{align*}
for all  $k \in \mathbb{Z}_{> 0}$, where $V_{k\varpi_{n-m}}$ is the $\OO_q(SU_n)$-comodule corresponding to the partition $(k, k, \dots, k)$ of length $n-m$. 
\end{thm}
We then use this theorem to give a novel  presentation of $S_q(\mathrm{Gr}_{n,m})$ generalising the classical ample bundle presentation of the homogeneous coordinate ring of the Grassmannian.

Finally, in \textsection 7 we summarise some recent applications of the Borel--Weil theorem to the study of the noncommutative K\"ahler structures of the quantum Grassmannians \cite{DOSFred,DOKSS,OSV}. In particular, we show how the general theory of K\"ahler structures allows us to extend the Borel--Weil theorem for positive line bundles to a quantum generalisation of the Bott--Borel--Weil theorem, and how to produce Dolbeault--Dirac Fredholm modules by tensoring with a negative line bundle.


\subsubsection*{Acknowledgements:} We would like to thank Edwin Beggs, Tomasz  Brzezi\'nski, Freddy D\'iaz, Matthias Fischmann, Andrey Krutov, Petr  Somberg, Karen Strung, Adam-Christiaan van Roosmalen, and Elmar Wagner, for useful discussions.


\section{Preliminaries on Differential Calculi and Quantum Principal Bundles}

We begin with a presentation of the necessary results on quantum group noncommutative geometry, namely the theory of covariant differential calculi and quantum principal bundles. All this material is by now quite well-know, and a more detailed presentation can be found in the recent monograph \cite{BeggsMajid:Leabh}.

\subsection{Differential Calculi and Complex Structures}

We begin a presentation of the general theory of differential calculi, complex structures, and their covariant versions over comodule algebras.

\subsubsection{Differential Calculi}

A {\em differential calculus} $\big(\Omega^\bullet \simeq \bigoplus_{k \in \bZ_{\geq 0}} \Omega^k, \exd\big)$ is a differential graded algebra (dg-algebra)  which is generated in degree $0$ as a dg-algebra, that is to say, it is generated as an algebra by the elements $a, \exd b$, for $a,b \in \Omega^0$. 
For a given algebra $B$, a differential calculus {\em over} $B$ is a differential calculus such that $B = \Omega^0$.
A {\em differential $*$-calculus} over a $*$-algebra $B$ is a differential calculus over $B$ such that the \mbox{$*$-map} of $B$ extends to a (necessarily unique) conjugate linear involutive map $*:\Omega^\bullet \to \Omega^\bullet$ satisfying $\exd(\omega^*) = (\exd \omega)^*$, and 
\begin{align*}
\big(\omega \wed \nu\big)^*  =  (-1)^{kl} \nu^* \wed \omega^*, &  & \text{ for all } \omega \in \Om^k, \, \nu \in \Omega^l. 
\end{align*}

A differential calculus  $\Omega^\bullet$ over a left $A$-comodule algebra $P$ is said to be {\em covariant} if the coaction $\DEL_L:P \to A \otimes P$ extends to a (necessarily unique) $A$-comodule algebra structure  $\DEL_L:\Om^\bullet \to A \otimes \Omega^\bullet$, with respect to which the differential $\exd$  is a left $A$-comodule map.  Covariance for a right $A$-comodule algebra is defined analogously. See \cite[\textsection 1]{BeggsMajid:Leabh} for a more detailed discussion of differential calculi.

\subsubsection{First-Order Differential Calculi}

A {\em first-order differential calculus} over an algebra $B$ is a pair $(\Om^1,\exd)$, where $\Omega^1$ is a $B$-bimodule and $\exd: B \to \Omega^1$ is a linear map for which the {\em Leibniz rule} holds
\begin{align*}
\exd(ab)=a(\exd b)+(\exd a)b,&  & a,b \in B,
\end{align*}
and for which $\Om^1$ is generated as a left $B$-module by those elements of the form~$\exd b$, for~$b \in B$. The {\em universal first-order differential calculus} over $B$ is the pair
$(\Om^1_u(B), \exd_u)$, where $\Om^1_u(B)$ is the kernel of the multiplication map $m_B: B \otimes B \to B$ endowed
with the obvious bimodule structure, and $\exd_u$ is the map defined by
\begin{align*}
\exd_u: B \to \Omega^1_u(B), & & b \mto 1 \otimes b - b \otimes 1.
\end{align*}
Every first-order differential calculus over $B$ is of the form   $\left(\Omega^1_u(B)/N, \,\proj \circ \exd_u\right)$, where $N$ is a $B$-sub-bimodule of $\Omega^1_u(B)$, and  
$$
\proj:\Omega^1_u(B) \to \Omega^1_u(B)/N
$$
is the canonical quotient map. This gives a bijective correspondence between calculi and sub-bimodules of $\Omega^1_u(B)$. For any subalgebra $B' \sseq B$, the \emph{restriction} of a first-order differential calculus over $B$ to $B'$ is the first-order differential calculus $\Omega^1(B') \sseq \Omega^1(B')$ over $B'$ generated by the elements $\mathrm{d}b$, for $b \in B'$.

We say that a differential calculus $(\Gamma^\bullet,\exd_{\Gamma})$ {\em extends} a first-order calculus $(\Omega^1,\exd_{\Omega})$ if there exists a bimodule isomorphism $\phi:\Omega^1 \to \Gamma^1$ such that  $\exd_{\Gamma} = \phi \circ \exd_{\Omega}$. Any first-order calculus admits an extension $\Omega^\bullet$ which is maximal in the sense that there exists a unique differential map from $\Omega^\bullet$ onto any other extension of $\Omega^1$, see  \cite[\textsection 1.5]{BeggsMajid:Leabh} for details. We call this extension the {\em maximal prolongation} of the first-order calculus.

\subsubsection{Complex Structures}

We now recall the definition of a complex structure as introduced in \cite{ BS,KLVSCP1}. This  abstracts the properties of the de Rham complex of a classical complex manifold \cite{HUY}. 

\begin{defn}\label{defnnccs}
A {\em complex structure} $\Om^{(\bullet,\bullet)}$ for a  differential $*$-calculus  $(\Om^{\bullet},\exd)$ is an $\mathbb{Z}_{\geq 0}^2$-algebra grading $\bigoplus_{(a,b)\in \mathbb{Z}_{\geq 0}^2} \Om^{(a,b)}$ for $\Om^{\bullet}$ such that, for all $(a,b) \in \mathbb{Z}_{\geq 0}^2$ it holds that 
\begin{enumerate}
\item \label{compt-grading}  $\Om^k = \bigoplus_{a+b = k} \Om^{(a,b)}$,
\item  \label{star-cond} $\big(\Om^{(a,b)}\big)^* = \Om^{(b,a)}$.
\item $\exd \Om^{(a,b)} \subseteq \Omega^{(a+1,b)} \oplus \Omega^{(a,b+1)}$,  for all $(a,b) \in \mathbb{Z}_{\geq 0}$.
\end{enumerate}
\end{defn}

We call an element of $\Om^{(a,b)}$ an $(a,b)$-form. For the  projections 
\begin{align*}
\proj_{\Om^{(a+1,b)}} : \Omega^{a+b+1} \to \Omega^{(a+1,b)}, & & \proj_{\Om^{(a,b+1)}} : \Omega^{a+b+1} \to \Omega^{(a,b+1)},
\end{align*}
we denote
\begin{align*}
\del|_{\Om^{(a,b)}} : = \proj_{\Om^{(a+1,b)}} \circ \exd, & & \ol{\del}|_{\Om^{(a,b)}} : = \proj_{\Om^{(a,b+1)}} \circ \exd.
\end{align*}
These definitions imply the identities 
\begin{align*}
\exd = \del + \adel, & &  \adel \circ \del = - \, \del \circ \adel, & & \del^2 = \adel^2 = 0,
\end{align*}
which is to say,  $\big(\!\bigoplus_{(a,b)\in \mathbb{Z}^2_{\geq 0}}\Om^{(a,b)}, \del,\ol{\del}\big)$ is a double complex.   Moreover,  both $\del$ and $\adel$ satisfy the graded Leibniz  rule and one has $\del(\omega^*) = \big(\adel \omega\big)^*$, and $\adel(\omega^*) = \big(\del \omega\big)^*$,  for all  $\omega \in \Omega^\bullet$. The {\em opposite} complex structure of a complex structure $\Om^{(\bullet,\bullet)}$ is the \mbox{$\mathbb{Z}^2_{\geq 0}$}-algebra grading  $\overline{\Om}^{(\bullet,\bullet)}$, defined by $\ol{\Om}^{(a,b)} := \Om^{(b,a)}$, for $(a,b) \in \mathbb{Z}_{\geq 0}^2$. 

For $\Om^\bullet$ a covariant differential $*$-calculus $\Om^\bullet$ over a left $A$-comodule algebra $P$,  we say that a complex structure for $\Om^\bullet$ is {\em covariant} if the $\mathbb{Z}_{\geq 0}^2$-decomposition is a decomposition in $^A\mathrm{Mod}$, which is to say, if $\Om^{(a,b)}$ is a left $A$-sub-comodule of $\Om^\bullet$, for each $(a,b) \in \mathbb{Z}_{\geq 0}$. 
A direct consequence of covariance is that the maps $\del$ and $\adel$ are left $A$-comodule maps.

\subsubsection{Connections and Holomorphic Structures}

We now recall the definition of a noncommutative holomorphic vector bundle, as defined in \cite{BS,KLVSCP1}.  A \emph{vector bundle} over $B$ will mean a finitely generated projective left $B$-modules. A \emph{line bundle} over $B$ will be an invertible $B$-bimodule $\EE$, where \emph{invertible} means that there exists another $B$-bimodule $\EE^{\vee}$ such that 
$
\EE \otimes_B \EE^{\vee} \simeq \EE^{\vee}\! \otimes_B \EE \simeq B.
$
Note that any such $\EE$ is automatically projective as a left $B$-module, and so, a line bundle is automatically a vector bundle. 
Building on this idea we define noncommutative holomorphic vector bundles via the classical Koszul--Malgrange characterisation of holomorphic bundles \cite{KM}. 

For $\Omega^\bullet$ a differential calculus over an algebra $B$, and $\mathcal{F}$ a left $B$-module, a \emph{connection} on $\F$ is a $\mathbb{C}$-linear map $\nabla:\mathcal{F} \to \Omega^1 \otimes_B \F$ satisfying 
\begin{align*}
\nabla(bf) = \exd b \otimes f + b \nabla f, & & \textrm{ for all } b \in B, f \in \F.
\end{align*}

With respect to a choice $\Omega^{(\bullet,\bullet)}$ of complex structure on $\Omega^{\bullet}$, a \emph{$(0,1)$-connection for $\mathcal{F}$} is a connection with respect to the differential calculus $(\Omega^{(0,\bullet)},\adel)$.

Any connection can be extended to a map $\nabla: \Omega^\bullet \otimes_B \mathcal{F} \to   \Omega^\bullet \otimes_B \mathcal{F}$ uniquely defined by 
\begin{align*}
\nabla(\omega \otimes f) =   \exd \omega \otimes f + (-1)^{|\omega|} \, \omega \wedge \nabla f, & & \textrm{for } f \in \F, \, \omega \in \Omega^{\bullet},
\end{align*}
for a homogeneous form $\omega$ with degree  $|\omega|$. The \emph{curvature} of a connection is the left $B$-module map $\nabla^2: \mathcal{F} \to \Omega^2 \otimes_B \mathcal{F}$. A connection is said to be {\em flat} if $\nabla^2 = 0$. Since $\nabla^2(\omega \otimes f) = \omega \wedge \nabla^2(f)$, a connection is flat if and only if  the pair $(\Omega^\bullet \otimes_B \F, \nabla)$ is a complex. 

\begin{defn}
For an algebra $B$, a \emph{holomorphic vector bundle over $B$} is a pair $(\mathcal{F},\adel_{\mathcal{F}})$, where  $\mathcal{F}$ is a finitely generated projective left $B$-module, and the map $\adel_{\mathcal{F}}: \mathcal{F} \to \Omega^{(0,1)} \otimes_B \mathcal{F}$ is a flat $(0,1)$-connection, which we call the \emph{holomorphic structure} for $(\F, \adel_{\F})$. 
\end{defn}

Note that for any fixed $a \in \mathbb{Z}_{\geq 0}$, a holomorphic vector bundle $(\F,\adel_{\F})$ has a naturally associated complex
$$
\adel_{\F}: \Omega^{(a,\bullet)} \otimes_B \F  \to \Omega^{(a,\bullet)} \otimes_B \F.
$$
For any $b \in \mathbb{Z}_{\geq 0}$, we denote by $H^{(a,b)}_{\adel}(\F)$ the $b^{\mathrm{th}}$-cohomology group of this complex.


\subsection{Principal Comodule Algebras} \label{subsection:PCA}

A right $H$-comodule algebra $(P,\Delta_R)$ is said to be a {\em $H$-Hopf--Galois extension of} $B := P^{\co(H)}$ 
if for  $m_P:P \otimes_B P \to P$ the multiplication of~$P$,
a bijection is given by 
$$
\can := (m_P \otimes \id) \circ (\id \otimes \DEL_R): P \otimes_B P \to P \otimes H.
$$
The Hopf--Galois condition is equivalent to exactness of the sequence
\begin{align} \label{eqn:qpbexactseq}
0 \longrightarrow P\Om^1_u(B)P {\buildrel \iota \over \longrightarrow} \Om^1_u(P) {\buildrel {\mathrm{ver}}\over \longrightarrow} P \oby H^+ \longrightarrow 0,
\end{align}
where $\Om^1_u(B)$ is the restriction of $\Om^1_u(P)$ to $B$, $\iota$ is the inclusion map, and 
$$
\mathrm{ver} := \mathrm{can} \circ \mathrm{proj}_B
$$
with $\proj_B$ the restriction to $\Omega^1_u(P)$ of the canonical projection $P \otimes P \to P \otimes_B P$. We call the map $\mathrm{ver}$ the \emph{vertical projection}.

We say that $P$ is faithfully flat as a right  $B$-module if the functor 
$$
P \otimes_B -: {}_B\mathrm{Mod} \to \mathrm{Mod_{\mathbb{C}}}
$$
preserves and reflects exact sequences. Faithful flatness as a left $B$-module is defined analogously.

\begin{defn}\
A {\em principal right $H$-comodule algebra}  is a right $H$-comodule algebra $(P,\DEL_R)$  such that  $P$ is a Hopf--Galois extension of $B := P^{\co(H)}$ and $P$ is faithfully flat as a right and left $B$-module.
\end{defn}

In most cases, directly verifying the requirements of this definition proves impractical. The following theorem gives a workable reformulation. 

\begin{defn}\label{defn:ell}
For a right $H$-comodule algebra $(P,\Delta_R)$, a {\em principal $\ell$-map} is a linear map $\ell:H \to P \otimes P$ satisfying 
\bet
\item $\ell(1_H) = 1_P \oby 1_P$,
\item $m_P \circ \ell = \e_H 1_P$,
\item $(\ell \oby \id_H) \circ \DEL_H = (\id_P \oby \DEL_R) \circ \ell$,
\item $(\id_H \oby \ell) \circ \DEL_H = (\DEL_L \oby \id_P) \circ \ell$,
\eet
where, for $\text{flip}: P \oby H \to H \oby P$  the flip map, $\DEL_L : = (S \oby \id) \circ \text{flip} \circ \Delta_R$.
\end{defn}
As shown in \cite{BrzBohm,BRzHajComptesT}, a right  $H$-comodule algebra $(P,\Delta_R)$ is principal if and only if  there exists a principal $\ell$-map $\ell: H \to P \oby P$.

We finish with a special type of comodule algebra which plays an important role in this paper. Let $\pi : A \to H$ be a surjective Hopf algebra map between Hopf algebras $A$ and $H$. The coaction
$
\Delta_R := (\id \otimes \pi) \circ \Delta : A \to A \otimes H,
$
gives $A$ the structure of a right $H$-comodule algebra. The associated space of coinvariant elements~$B:=A^{\co(H)}$ is called a \emph{quantum homogeneous space} if $A$ is faithfully flat as a right $B$-module. As is well-known, if $H$ is comsemisimple then faithful flatness is automatic, in fact the pair $(A,\Delta_{R})$ is a principal comodule algebra (see \cite[\textsection 3.3]{DOKSS} for a more detailed discussion). Note that for the special case of $\pi=\e:A \to \mathbb{C}$, which we call the \emph{trivial quantum homogeneous space}, it holds that $B = A$. 

\subsection{Quantum Principal Bundles}\label{subsection:QPB}

 A non-universal first-order differential calculus on~$P$ is said to be {\em right $H$-covariant} if the following (necessarily unique) map is well defined
\begin{align*}
\DEL_R: \Om^1(P) \to \Om^1(P) \otimes H, & & p \exd q \mto p_{(0)} \exd q_{(0)} \otimes p_{(1)} q_{(1)}.
\end{align*} 
(Covariance for a first-order differential calculus over a left comodule algebra is defined similarly.) The following definition, due to Brzezi\'nski and Majid  \cite{BeggsMajid:Leabh,TBSM1}, presents sufficient criteria for the existence of a non-universal version of the sequence \eqref{eqn:qpbexactseq}.

\begin{defn} \label{qpb}
Let $H$ be a Hopf algebra. A {\em  quantum principal $H$-bundle} is a pair $(P,\Omega^1(P))$, consisting of a right $H$-comodule algebra $(P,\Delta_R)$ and a right-$H$-covariant calculus $\Om^1(P)$, such that:
\begin{enumerate}
    \item $P$ is a Hopf--Galois extension of $B  = P^{\,\co(H)}.$
    \item If $N \sseq \Om^1_u(P)$ is the sub-bimodule of the universal calculus corresponding to $\Om^1(P)$, we have $\mathrm{ver}(N) = P \otimes I$, for some $\mathrm{Ad}$-sub-comodule right ideal 
\[ I \sseq H^+ :=  \ker(\e: H \to \mathbb{C}),
\] where $
\mathrm{Ad} : H \to H \otimes H$ is defined by $\mathrm{Ad}(h) := h_{(2)} \otimes S(h_{(1)}) h_{(3)}.$
\end{enumerate}
\end{defn}

Denoting by  $\Om^1(B)$ the restriction of $\Om^1(P)$ to $B$, and $\Lambda^1(H) := H^+/I$, the quantum principal bundle definition implies that an exact sequence is given by
\begin{align} \label{Eqn:qpbexactseq}
0 \longrightarrow P\Om^1(B)P {\buildrel \iota \over \longrightarrow} \, \Om^1(P) {\buildrel {~\mathrm{ver}~~}\over \longrightarrow} P \otimes \Lambda^1(H) \longrightarrow 0,
\end{align}
where by abuse of notation $\mathrm{ver}$ denotes the map induced on $\Om^1(P)$ by identifying~$\Om^1(P)$ as a quotient of $\Om_u^1(P)$.

\begin{defn}
A {\em principal connection} for a quantum principal $H$-bundle $(P,\Omega^1(P))$ is a left $P$-module, right $H$-comodule, projection $\Pi:\Om^1(P) \to \Om^1(P)$ satisfying
$$
\ker(\Pi) = P\Om^1(B)P.
$$
A principal connection $\Pi$ is called {\em strong} if $(\id - \Pi) \big(\exd P\big) \sseq \Om^1(B)P$. 
\end{defn}
The existence of a principal connection is equivalent to the existence of a left \linebreak $P$-module, right $H$-comodule, splitting $s:P \otimes H^+ \to \Omega^1(P)$ of the vertical projection $\mathrm{ver}$. Explicitly, the equivalence is determined by 
\begin{align} \label{eqn:connectionsplittlin}
\Pi(\omega) = s(\mathrm{ver}(\omega)), & & \textrm{ for } \omega \in \Omega^1(P).
\end{align}
Moreover, the existence of a strong principal connection for a comodule algebra is equivalent to the comodule algebra being principal \cite{BrzBohm},\cite{BRzHajComptesT}. For a Hopf--Galois extension $B = P^{\co(H)}$,  a principal $\ell$-map  $\ell: H \to P \oby P$, a strong principal connection  is defined by
\begin{align} \label{elltopi}
\Pi_{\ell}:= (m_P \otimes \id) \circ (\id \otimes \ell) \circ \mathrm{ver}: \Omega^1_u(P) \to \Omega^1_u(P).
\end{align}
Moreover, this association gives a bijective correspondence between principal $\ell$-maps  and strong principal connections.
(For an expository presentation of this equivalence, as well as the general theory of principal comodule algebras, we recommend the notes \cite[\textsection VII]{TBGS}.) 


\subsection{Connections from Principal Connections} \label{subsection:ConnFromPrinCs}

Denoting by ${}_P\mathrm{Mod}$ the category of left $P$-modules, we define a functor
\begin{align*}
\Psi:{}^{H}\mathrm{Mod} \to {}_P\mathrm{Mod}, & & V \mapsto P \square_H V,
\end{align*}
which acts on morphisms $\gamma:V \to W$ as $\Psi(\gamma) = \id \otimes \gamma$. For any $\F := \Psi(V)$, we can use principal connections to define a connection $\nabla:\F \to \Omega^1(B) \otimes_B \F$. 
Note first that we have a natural embedding 
\begin{align*}
j:\Om^1(B) \oby_B \F \hookrightarrow  \Om^1(B)P \, \square_H V,
\end{align*}
given by the multiplication map.
A strong principal connection $\Pi$ defines a connection $\nabla$ on $\F$ by
\begin{align*}
 \nabla: \F \to  \Om^1(B) \oby_B \F, &   & \sum_i p_i \otimes v_i \mto  j^{-1} \big((\id - \Pi)(\exd p_i) \oby v_i \big).
\end{align*}
Indeed, since $\exd$ and the projection $\Pi$ are both right \mbox{$H$-comodule} maps, their composition $(\id - \Pi) \circ \mathrm{d}$ is a right $H$-comodule map. Hence the image of $(\id - \Pi) \circ \mathrm{d}$  is contained in  $j\left(\Om^1(B) \oby_B \F\right)$ and $\nabla$ defines a connection. 

Finally, we consider the special case where $P$ is endowed with a left $A$-coaction giving it the structure of an $(A,H)$-bicomodule. If the principal connection $\Pi$ is a left $A$-comodule map, then we see that the connection $\nabla$ will also be left $A$-comodule map. This will be the case for all the connections considered in this paper. 

\subsection{Some Categorical Equivalences}

We denote by~${}^A_B\mathrm{Mod}$ the category of \emph{relative Hopf modules}, that is, the category whose  objects \mbox{$\DEL_L:\mathcal{F} \to A \otimes \mathcal{F}$} are left \mbox{$A$-comodules},  endowed with a left $B$-module structure such that,  for all $f \in \mathcal{F},$ \mbox{$b \in B$}, we have $\DEL_L(bf) = \Delta_L(b)\DEL_L(f)$, 
and whose morphism are left $A$-comodule, left $B$-module, maps. 
Moreover, we denote by ${}^H\mathrm{Mod}$  the category whose objects are left \mbox{$H$-comodules}, and whose morphisms are left $H$-comodule maps. 

Denoting $B^+ := B \cap \ker(\e)$, consider the functor
$
\Phi:{}^A_B\mathrm{Mod} \to {}^H\mathrm{Mod},
$
defined by $\Phi(\mathcal{F}) := \mathcal{F}/B^+\mathcal{F}$,
where the left $H$-comodule structure of $\Phi(\mathcal{\F})$ is given by 
\begin{align}\label{eqn:TakleftHcoaction}
\Delta_L[f] := \pi(f_{(-1)})\otimes [f_{(0)}]
\end{align}
(with square brackets denoting the  coset of an element  in $\Phi(\mathcal{\F})$.) In the other direction, using the cotensor product we can define the functor $\Psi:  {}^H\mathrm{Mod} \to {}^A_B\mathrm{Mod}$ by setting $\Psi(V) := A \,\square_H V$,  where the left $B$-module and left $A$-comodule structures of $\Psi(V)$ are defined on the first tensor factor, and if $\gamma$ is a morphism in ${}^H\mathrm{Mod}$, then $\Psi(\gamma) := \id \otimes \gamma$.

As established in~\cite[Theorem 1]{Tak}, an adjoint equivalence of categories between~${}^A_B\mathrm{Mod}$ and~${}^H\mathrm{Mod}$, which we call \emph{Takeuchi's equivalence}, is given by the functors $\Phi$ and $\Psi$ with the unit natural isomorphism
\begin{align} \label{eqn:unit}
\unit: \F \to \Psi \circ \Phi(\F), & & f \mto f_{(-1)} \otimes [f_{(0)}],
\end{align}
and the counit natural isomorphism
\begin{align} \label{eqn:counit}
\qquad \counit:\Phi \circ  \Psi(V) \to V, \,  & & \Big[\sum_i a^i \otimes v^i\Big] \mto \sum_{i} \e(a^i)v^i.
\end{align} 
In what follows, when it is not clear from the context, we will put a subscript to Takeuchi's functors to specify in which module category we are working, so for example $\Phi_P$ will denote Takeuchi's functor $\Phi$ with domain $^A_P\mathrm{Mod}.$\\
Consider next the category ${}^A_B\mathrm{Mod}_B$ whose objects are formed by endowing objects $\mathcal{F} \in {}^A_B\mathrm{Mod}$ with a right $B$-module structure, giving $\F$ the structure of a $B$-bimodule satisfying $\Delta_L(fb) = \Delta_L(f)\Delta(b)$, for all $b \in B$, $f \in \mathcal{F}$, and whose morphisms are left $A$-comodule and $B$-bimodule maps. Consider also the category ${}^H\mathrm{Mod}_B$ whose objects are left $H$-comodules $\Delta_L:V \to H \otimes V$ and right $B$-modules, satisfying $\Delta_L(vb) = v_{(-1)}\pi(b_{(1)}) \otimes v_{(0)}b_{(2)}$, and whose morphisms are left $H$-comodule and right $B$-module maps. For any $\F \in {}^A_B\mathrm{Mod}_B$, we endow $\Phi(\F)$ with the right $B$-module structure $[f]b := [fb]$, for $f \in \mathcal{F}$, $b \in B$. Moreover, we endow $\Psi(V)$ with a right $B$-module structure by setting
\begin{align*}
\left(\sum_i a_i \otimes v_i\right)\!b := \sum_i a_ib_{(1)} \otimes \left(v_i b_{(2)}\right)\!, 
\end{align*}
for any $b \in B$, and any $\sum_i a_i \otimes v_i  \in  A \square_H V$.
These definitions allow us to extend Takeuchi's equivalence to an equivalence between ${}^A_B\mathrm{Mod}_B$ and ${}^H\mathrm{Mod}_B$. 

Denote by  $^A_B\textrm{Mod}_0$ the full sub-category of $^A_B\mathrm{Mod}_B$ whose objects $\F$  satisfy the identity 
$
\mathcal{F}B^+ = B^+\mathcal{F}.
$ 
Moreover, denote by ${}^H\mathrm{Mod}_0$ the full subcategory of ${}^H\mathrm{Mod}_B$ consisting of those objects endowed with the trivial right $B$-action, which is to say, those objects $V$ for which $v \tl b = \e(b)v$, for all $v \in V$, and $b \in B$. (This category is clearly isomorphic to $^H\mathrm{Mod}$.) As explained in \cite[Corollary 2.5]{MMF3}, Takeuchi's equivalence restricts to an equivalence between $^A_B\textrm{Mod}_0$ and $^H\mathrm{Mod}$. 

For the special case of the trivial quantum homogeneous space we get an equivalence between the categories ${}^A_A\mathrm{Mod}_A$ and $\mathrm{Mod}_A$. In this case we denote $\Phi$ by $F$, and $\Psi$ by $A \otimes -$ (note that for this special case the cotensor product reduces to the usual tensor product over $\mathbb{C}$.)  When discussing bicovariant differential calculi in \textsection \ref{subsection:bicovariant} we find it useful to consider the category  of \emph{Hopf bimodules} ${}^A_A\mathrm{Mod}^A_A$, whose objects are objects of ${}^A_A\mathrm{Mod}_A$ endowed with a right $A$-comodule structure satisfying the obvious compatibility relations, and whose morphims are $A$-bimodule, $A$-bicomodule, maps. This category is equivalent to $\mathcal{YD}^A_A$ the category of  \emph{Yetter--Drinfeld modules} over $A$, which is to right $A$-module right $A$-comodules satisfying 
$$
v_{(0)}a_{(1)} \otimes v_{(1)}a_{(2)} = (vh_{(2)})_{(0)} \otimes h_{(1)}(vh_{(2)})_{(1)} 
$$
together with right $A$-module right $A$-comodule morphisms. Explicitly, for any $\mathcal{F} \in {}^A_A\mathrm{Mod}^A_A$ we endow $F(\mathcal{F})$ with the right \emph{adjoint} $A$-comodule structure
\begin{align}
F(\mathcal{F}) \to F(\mathcal{F}) \otimes A, & &  [f] \mapsto [f_{(0)}] \otimes S(f_{(-1)})f_{(1)},
\end{align}
and for $V \in \mathcal{YD}^A_A$, we endow $A \otimes V$ with the tensor product right $A$-coaction.

\subsection{Covariant First-Order Differential Calculi} \label{subsection:CovariantFODCi}

For $A$ a Hopf algebra, and $B$ a left $A$-comodule algebra, a first-order differential calculus $\Omega^1(B)$ over $B$ is said to be {\em left $A$-covariant} if there exists a left \mbox{$A$-coaction $\DEL_L: \Omega^1(B) \to A \otimes \Omega^1(B)$} (necessarily unique) satisfying
\begin{align} \label{eqn:FODCLeftCov}
\DEL_L(b\exd b') = \DEL_L(b) (\id \otimes \exd)\DEL_L(b'),  &  & \text{ for all } b,b' \in B.
\end{align}
A differential calculus over $B$ is said to be {\em left $A$-covariant} if it admits a left $A$-coaction satisfying the obvious analogue of \eqref{eqn:FODCLeftCov}. The maximal prolongation of a covariant first-order calculus is automatically covariant.

For any quantum homogeneous space $B = A^{\co(H)}$, a covariant first-order differential calculus over $B$ is naturally an object in $^A_B\mathrm{Mod}_B$.
Moreover, an isomorphism in ${}^H\mathrm{Mod}_B$ is given by 
\begin{align}\label{eqn:basicsigma}
\sigma:\Phi(\Omega_u^1(B)) \to B^+, & & [b\exd b'] \mapsto \e(b)(b')^+.
\end{align}
As established in \cite{Herm} (see also \cite{Maj}) we can use $\sigma$ to classify covariant first-order differential calculi in terms of subobjects $I \sseq B^+$ in ${}^H\mathrm{Mod}_B$. Explicitly, for a $B$-subbimodule $N \sseq \Omega^1_u(B)$ the corresponding ideal is given by 
$
\sigma(\Phi(N)) = I.
$
This gives us 
the commutative diagram  
\begin{align*}
\xymatrix{ 
\Omega^1(B)     \ar[rrr]^{(\id \otimes \sigma) \circ \unit}                  & & &    A \square_H V(B),  \\
B   \,\, \ar^{\exd}[u]  \ar[rrru]_{(\id  \otimes [\cdot]) \circ \Delta_L}                     & & &
}
\end{align*}
where $V(B)=B^+/I$ and $[\cdot]: B^+ \to V(B) $ denotes the canonical projection map. Moreover by abuse of notation $\sigma$ denotes the projection of the map \eqref{eqn:basicsigma} to the non-universal calculus $\Omega^1(B)$. 

For the special case $B=A$, 
the classification of covariant first-order differential calculi discussed above reduces to Woronowicz's celebrated classification \cite{WoroDC} which gives a bijection between covariant first-order differential calculi over $A$ and right ideals of $A^+$.

Consider next the restriction of a left $A$-covariant first-order differential calculus $\Omega^1(A)$ to a  (necessarily left $A$-covariant) first-order differential calculus on $B$. This gives us the following commutative diagram in the category ${}^A\mathrm{Mod}$:
\begin{align*}
\xymatrix{ 
\Omega^1(A)  \ar[rrr]^{(\id \otimes \sigma) \, \circ \, \unit_{\Omega^1(A)} }                    & & &  A \otimes \Lambda^1(A)   \\
\Omega^1(B) \ar[u]^{\iota}  \ar[rrr]_{(\id \otimes \sigma) \, \circ \, \unit_{\Omega^1(B)}} \,\,    & & &   \ar[u]_{\id \otimes \iota'} A \square_H V(B),
}
\end{align*} 
where $\iota:\Omega^1(B) \to \Omega^1(A)$ and $\iota':V(B) \to \Lambda^1(A)$ are the evident inclusions. Recall that the antipode of $H$ implements an equivalence between the
categories of left and right $H$-comodules. Explicitly, any $V \in {}^H\mathrm{Mod}$ can be endowed with a right $H$-coaction $\Delta_{R,H}$ defined by 
\begin{align}
    \Delta_{R,H}(v) := v_{(0)} \otimes S(v_{(-1)}), & & \textrm{ for any } v \in V.
\end{align}
An important point to note is that endowing $V(B)$ with its right $H$-comodule structure makes the inclusion 
$$
\iota': V(B) \to \Lambda^1(A)
$$
a right $H$-comodule map. 

\subsection{Coquasi-Triangular Structures and Bicovariant Calculi} \label{subsection:bicovariant}

Continuing with the special case of $B=A$, we say that a first-order differential calculus $\Omega^1(A)$ is {\em bicovariant} if it is both left and right $A$-covariant and satisfies
$
(\id \otimes \DEL_{R,A}) \circ \DEL_{L,A} = (\DEL_{L,A} \otimes \id) \circ \DEL_{R,A},
$ 
Using the equivalence of the category of Hopf bimodules and the category of Yetter--Drinfeld modules it is possible to identify those ideals $I \sseq A^+$ that give bicovariant calculi \cite[Theorem 1.8]{WoroDC}. Explicitly, a left-covariant first-order calculus $\Omega^1(A)$ is bicovariant if and only if the corresponding ideal $I \sseq A^+$ is a  subcomodule of $A^+$ with respect to the {\em (right) adjoint coaction} 
\begin{align*}
\mathrm{Ad_R}: A^+ \to A^+ \otimes A, & & a \mapsto a_{(2)} \otimes S(a_{(1)})a_{(3)}.
\end{align*}
When considering left covariant first-order calculi over Hopf algebras we use the notation 
\begin{align} \label{asdf}
\Lambda^1(A) := A^+/I.
\end{align}

\begin{defn} \label{defn:CQT}
We say that a Hopf algebra $A$ is {\em coquasi-triangular} if it is equipped with a convolution-invertible linear map $r :A \otimes A \rightarrow \bC$ obeying, for all $a,b,c \in A$, the relations
\begin{align*} 
r(ab\otimes c) = r(a \otimes c_{(1)})r(b\otimes c_{(2)}),            &   & r(a\otimes bc) = r(a_{(1)}\otimes c)r(a_{(2)}\otimes b),
\end{align*}
\begin{align*} 
 r(a_{(1)}\otimes b_{(1)})a_{(2)} b_{(2)} =  r(a_{(2)}\otimes b_{(2)}) b_{(1)}a_{(1)}, & & r(a \otimes 1) = r(1 \otimes a) = \e(a).
\end{align*}
\end{defn}


For any $a \in A$, we have a linear map 
\begin{align*}
Q_a: A \to \mathbb{C}, & &  b  \mto  r(a_{(1)} \otimes b_{(1)})r(b_{(2)} \otimes a_{(2)}).
\end{align*}
For an $A$-comodule $V$, with coordinate space $C(V)$, consider the subspace
\begin{align*}
I_V := \left\{b \in A^+ \, | \, Q_a(b) = 0, \textrm{ for all } a \in C(V) \right\}\!.
\end{align*}
A direct calculation confirms that  $I_A$ is a right ideal of $A^+$. Moreover, $I_V$ is an $\mathrm{Ad}$-subcomodule of $A^+$ for every $A-$comodule $V$ \cite[\textsection 10.1.3]{KSLeabh} . Hence there exists an associated bicovariant calculus $\Om^1_{V}(A)$ associated to every $V$. 




\section{Principal Pairs and Takeuchi's Equivalence}

In this section we introduce the notion of a principal pair of Hopf algebras, and show how one can construct quantum principal bundles from them. This serves as our formal framework for investigating the quantum Grassmannian sphere the quantum Grassmannians in \textsection \ref{section:HKRestriction} and \textsection \ref{section:BorelWeil}.

\subsection{Principal Pairs}

In this subsection we introduce the notion of a principal pair of surjective Hopf algebras maps. As explained below, this generalises the definition of a homogeneous principal comodule algebra, and comes with an associated generalisation of Takeuchi's equivalence. It provides the general framework in terms of which the results of \textsection \ref{Subsection:PPApplications} below are presented. We begin with a technical lemma necessary for the statement of the definition of a principal pair. 

\begin{lem} \label{lem:PPsetup}
Let $A$ and $K$ be two Hopf algebras, and $\rho:A \to K$ a Hopf algebra map, with associated right $K$-coaction $\Delta_{R,K}: A \to A \otimes K$. For any 
right $K$-subcomodule subalgebra $P \sseq A$, it holds that 
\begin{enumerate}
    \item $T := \rho(P)$ is a sub-bialgebra of $K$,
    \item $T$ is the smallest subspace of $K$ for which $\Delta_{R,K}(P) \sseq P \otimes T$,
    \item in particular, $\Delta_{R,K}$ restricts to a coaction $\Delta_{R,T}:P \to P \otimes T$.
\end{enumerate}
\end{lem}
\begin{proof}
~~\\
1. ~~ Since $\rho$ is an algebra map, it is clear that $T$ is a subalgebra of $K$. Since $P$ is a left coideal of $A$, as well as a right $K$-subcomodule of $A$, we can write $\Delta(p)$, for any $p \in P$, as a linear combination of elements in $P \otimes P$ and elements in $A \otimes \ker(\rho)$. Thus for any $t = \rho(p) \in T$, we see that 
$$
\Delta(t) = \Delta(\rho(p)) = \rho(p_{(1)}) \otimes \rho(p_{(2)}) \in \rho(P) \otimes \rho(P) = T \otimes T.
$$

2. ~~ Since $P$ is a left coideal of $A$, we see that, for any $p \in P$,
$$
\Delta_{R,K}(p) = p_{(1)} \otimes \rho(p_{(2)}) \in P \otimes T.
$$
Thus $\Delta_{R,K}(P) \sseq P \otimes T$. Let $T' \sseq T$ be a subspace satisfying $\Delta_{R,K}(P) \sseq P \otimes T'$. For any $p \in P$, we see that
\begin{align*}
    \rho(p) = (\e \otimes \id)(p_{(1)} \otimes \rho(p_{(2)})) = (\e \otimes \id) \circ \Delta_{R,K}(p) \in T'.
\end{align*}
Thus $T \sseq T'$, meaning that $T' = T$ is the smallest subspace satisfying $\Delta_{R,K}(P) \sseq P \otimes T$. 

3. ~~ This follows directly from 2.
\end{proof}


\begin{defn}\label{defn:principapair}
For $A,H,K$ a triple of Hopf algebras, a \emph{principal pair} is a pair of surjective Hopf algebra maps $(\pi:A \to H,\rho:A \to K)$ such that
\begin{enumerate}
    \item $P:=A^{\co(H)}$ is a quantum homogeneous $A$-space, 
    \item $P$ is a right $K$-subcomodule of $A$,
    \item $T =\rho(P)$ is a Hopf subalgebra of $K$,
    \item denoting by $\Delta_{R,T}$ the restriction of $\Delta_{R,K}$ to $P$, the pair $(P,\Delta_{R,T})$ is a principal comodule algebra.
\end{enumerate}
\end{defn}

Note that for the special case of $H=K$ and $\pi=\rho$, we recover a homogeneous principal comodule algebra.

We now introduce a novel variation on Takeuchi's equivalence for principal pairs. 
Generalising the category of relative Hopf modules we have the following category.

\begin{defn}
Let $A,H,K$ be a triple of Hopf algebras, $(\pi:A \to H,\rho:A \to K)$ a \emph{principal pair} of Hopf algebra maps, and denote $P=A^{\co(H)}$. The objects of the category ${}^A_P\textrm{Mod}^K$ are objects  $\F \in {}^A_P\mathrm{Mod}$ endowed with a right $K$-comodule structure $\Delta_{R,K}:\F \to \F \otimes K$ giving $\F$ the structure of an $A$-$K$-bicomodule and satisfying 
$$
\Delta_{R,K}(pf) = \Delta(p)\Delta_{R,K}(f).
$$
The morphisms of the category are those morphisms in ${}^A_P\mathrm{Mod}$ which are right $K$-comodule maps.
\end{defn}

Generalising the category of left $H$-comodules we have the following category.

\begin{defn}
The objects of the category ${}^H\mathrm{Mod}^{K,\square}$ are left $H$-comodules, 
which are also right $K$-comodules (but not necessarily $H$-$K$-bicomodules) 
such that the cotensor product $A \square_H V$ is a right $K$-subcomodule of $A \otimes V$ equipped with the tensor product right $K$-coaction 
$$
\Delta_{\rho} \otimes \Delta_{R}: A \otimes V \to A \otimes V \otimes K,
$$
where $\Delta_{R,K}$ is the right $K$-coaction of $V$. The morphisms of the category are left $H$-comodule maps, which are also right $K$-comodule maps.
\end{defn}

We now introduce a functor 
$$
\Psi:{}^H\mathrm{Mod}^{K,\square} \to {}^A_P\textrm{Mod}^K
$$
by endowing $A \square_H V$ with the tensor product  right $H$-comodule, and by operating on a morphisms as for Takeuchi's functor $\Psi$.

\begin{prop} ~~~
\begin{enumerate}
\item  A functor $\Phi:{}^A_P\mathrm{Mod}^K \to {}^H\mathrm{Mod}^{K,\square}$ by considering $F \in {}^A_P\mathrm{Mod}^K$ as an object in ${}^A_P\mathrm{Mod}$  and then endowing $\Phi(F)$ with the right $K$-coaction defined by 
\begin{align*}
\Phi(\F) \to \Phi(\F) \otimes K, & & [f] \mapsto [f_{(0)}] \otimes  \rho\left(S(f_{(-1)})\right)f_{(1)},
\end{align*}
and by operating on morphisms as for Takeuchi's functor $\Phi$.

\item Takeuchi's unit and counit maps given in \eqref{eqn:unit} and \eqref{eqn:counit} are here right $K$-comodule maps, and hence a unit-counit equivalence between ${}^A_P\mathrm{Mod}^K$ and ${}^H\mathrm{Mod}^{K,\square}$ is given by $\Psi$, $\Phi$, $\unit$, and $\counit$.
\end{enumerate}
\end{prop}

\begin{proof}
~~~ \\
\begin{enumerate}
\item We need to show that the right $K$-action is well-defined. Take an element of the form $pf$ and note that 
$$
\Delta_{R,K}(pf) = p_{(2)}f_{(0)} \otimes \rho\left(S\left(p_{(1)}f_{(-1)}\right)\right)f_{(1)},
$$
where $\Delta_{R,K}: \F \to \F \otimes K$ is the obvious lift of $\Delta_{R,K}$ to a right $K$-coaction on $\F$. Since $A^+$ is closed under the right adjoint coaction, $P^+$ is closed under the right adjoint coaction, meaning that 
$$
[p_{(2)}f_{(0)}] \otimes \rho\left(S\left(p_{(1)}f_{(-1)}\right)\right)f_{(1)} = 0.
$$
Hence the right $K$-coaction is well-defined. To see that $\Phi(\F)$ is well-defined as an object in ${}^H\mathrm{Mod}^{K,\square}$ it remains to show that $A \square_H V$ is a right $K$-subcomodule of $A \otimes \Phi(\F)$. That this is true follows from the fact that, for any $f \in \mathcal{F}$, we have 
\begin{equation} \label{ale}\begin{split}
\Delta_{R,K}\!\left(f_{(-1)} \otimes [f_{(0)}]\right) = & \, f_{(-3)} \otimes [f_{(0)}] \otimes \rho(f_{(-2)}S(f_{(-1)}))f_{(1)} \\ 
= & \, f_{(-1)} \otimes [f_{(0)}] \otimes f_{(1)}\\ 
= & \, \unit(f_{(0)}) \otimes f_{(1)}.
\end{split}
\end{equation}

Note next that, for $\mathcal{F},\mathcal{G} \in {}^A_P\mathrm{Mod}^K$, and $f:\mathcal{F} \to \mathcal{G}$ a morphism, the fact that $f$ is a bicomodule map means that $\Phi(f)$ is a left $H$-comodule map, and a right $K$-comodule map, which is to say a morphism in ${}^H\mathrm{Mod}^{K,\square}$. Thus $\Phi$ is a well-defined functor. 

\item Equation \eqref{ale} implies that the following diagram commutes
\begin{align*}
\xymatrix{ 
f   \ar[d]_{\Delta_{R,K}}       \ar[rrr]^{\unit}                  & & &     f_{(-1)} \otimes [f_{(0)}]\ar[d]^{\Delta_{R,K}}  \\
f_{(0)} \otimes f_{(-1)}   \,\,   \ar[rrr]_{\unit \otimes \id}                     & & & f_{(0)} \otimes f_{(-1)} \otimes f_{(1)}.
}
\end{align*}
Thus we see that $\unit$ is a right $K$-comodule map. 

For $V \in {}^H\mathrm{Mod}^K$, $v \in V$, and $\sum_i a_i \otimes v_i$, we see that 
\begin{align*}
\Delta_{R,K}\!\left(\left[\sum_i a_i \otimes v_i\right]\right) = & \left[\sum_i (a_i)_{(2)} \otimes (v_i)_{(0)}\right] \otimes S((a_i)_{(1)})(a_i)_{(3)}v_{(1)}\\
= & \left[\sum_i a_i \otimes (v_i)_{(0)}\right]\! \otimes (v_i)_{(1)}. 
\end{align*}
This implies that the following diagram commutes
\begin{align*}
\xymatrix{ 
\left[\sum_i a_i \otimes v_i\right]    \ar[rrr]^{\counit}       \ar[d]_{\Delta_{R,K}}                  & & &    \sum_i \e(a_i)v_i \ar[d]^{\Delta_{R,K}}   \\
\left[\sum_i a_i \otimes (v_i)_{(0)}\right] \otimes (v_i)_{(1)}   \,\,   \ar[rrr]_{\counit \otimes \id}                     & & & \sum_i \e(a_i)(v_i)_{(0)} \otimes (v_i)_{(1)}.
}
\end{align*}
Hence $\counit$ is a right $K$-comodule map as claimed.
\end{enumerate}
\end{proof}

\subsection{The First-Order Differential Calculus Induced on $\Omega^1(B)$}

   In this subsection, for a principal pair $(\pi:A \to H,\rho:A\to K)$, we consider the restriction of a left $A$-covariant right $K$-covariant first-order differential calculus $\Omega^1(P)$ over $P$ to a  first-order differential calculus $\Omega^1(B)$ over $B=P^{\co(K)}$. We note that $\Omega^1(B)$ is necessarily left $A$-covariant, and that the right $K$-coaction of $\Omega^1(P)$ restricts to the trivial right $K$-coaction on $\Omega^1(B)$. We begin with the following commutative diagram generalising the special case presented in \textsection \ref{subsection:CovariantFODCi}. (The proof is an evident generalisation of the special case and hence is omitted.)

\begin{prop}\label{prop:commVPVB}
A commutative diagram in the category ${}^A\mathrm{Mod}$ is given by
\begin{align*}
\xymatrix{ 
\Omega^1(P)  \ar[rrr]^{(\id \otimes \sigma) \, \circ \, \unit_{\Omega^1(P)} }                    & & &  A \otimes V(P)   \\
\Omega^1(B) \ar[u]^{\iota}  \ar[rrr]_{(\id \otimes \sigma) \, \circ \, \unit_{\Omega^1(B)}} \,\,    & & &   \ar[u]_{\id \otimes \iota'} A \square_H V(B),
}
\end{align*} 
where $\iota:\Omega^1(B) \to \Omega^1(P)$ and $\iota':V(B) \to V(P)$ are the evident inclusions. Moreover, endowing $V(B)$ with its right $H$-comodule structure given by the antipode gives the inclusion 
$$
\iota': V(B) \to V(P)
$$
the structure of a right $K$-comodule map. 
\end{prop}

\begin{cor}
The subobject $I \sseq B^+$ in the category ${}^K\mathrm{Mod}_B$,  corresponding to the calculus $\Omega^1(B)$, is given by 
$$
I \sseq B^+ \cap \sigma(\Phi(N)),
$$
where $N$ is the sub-bimodule of $\Omega^1_u(P)$ corresponding to the calculus $\Omega^1(P)$.
\end{cor}
\begin{proof}
The identity follows from the fact that  
\begin{align*}
    \sigma\!\left(\Phi(N \cap \Omega_u^1(B))\right) = \sigma(\Phi(N)) \cap \sigma\!\left(\Phi(\Omega^1_u(B)\right),
\end{align*}
where in the second identity we have used the fact that $\Phi$ is an equivalence and hence commutes with pullbacks, and in particular intersections. 
\end{proof}

An important point to note is that if we start with a covariant calculus on $A$, induce a covariant calculus on $P$ by restriction, and then in turn restrict to a calculus on $B$, we can embed $\iota'(V(B)) \subseteq  V(P)$ the nested pair of subsequences into $\Lambda^1(A)$. This fact will be used tacitly throughout \textsection \ref{section:RestPresHK} when we present the Heckenberger--Kolb calculus in quantum principal bundle terms.

\subsection{Principal Pairs and Universal Quantum Principal Bundles}\label{Subsection:PPApplications}

Throughout this subsection, so as to avoid tedious repetition,  $A,H$, and $K$ will denote Hopf algebras, and $(\pi':A \to H,\rho:A\to K)$ a principal pair. The space of coinvariants is denoted by $P:=A^{\co(H)}$, and the equivalent pair of categories by ${}^A_P\mathrm{Mod}^K$ and ${}^H\mathrm{Mod}^{K,\square}$. Moreover, we write $T = \rho(P)$ and $B:=P^{\co(T)}$.

We start with a technical lemma establishing some necessary results about universal quantum principal bundles.

\begin{lem} \label{lem:PPver}
~~~
\begin{enumerate}
    \item The map $\mathrm{ver}$ is a morphism in the category ${}^A_P\mathrm{Mod}^K$.
    \item Endow $P^+$ with the left $H$-coaction $(\pi' \otimes \id) \circ \Delta$, and with the right adjoint $K$-coaction 
    $$
    (\id \otimes \rho) \circ \Ad_R: P^+ \to P^+ \otimes K.
    $$
    A left $H$-comodule, right $K$-comodule, isomorphism is given by 
    \begin{align*}
        \sigma:\Phi\!\left(\Omega^1_u(P)\right)
        \to P^+, & & [p\exd q] \mapsto \e(p)q^+.
         \end{align*}
    Hence $P^+$ is an object in ${}^H\mathrm{Mod}^{K,\square}$.
    \item Endowing $T^+$ with the trivial left $H$-coaction and the right adjoint coaction 
    $$
    \Ad_R:T^+ \to T^+ \otimes T \sseq T^+ \otimes K,
    $$
    it holds that $A \, \square_{H} T^+  = P \otimes T^+$ and hence $T^+$ is an object in ${}^H\mathrm{mod}^{K,\square}$. 
    
\item The map $\rho|_{P^+}:P^+ \to T^+$ is a morphism in the category ${}^H\mathrm{Mod}^{K,\square}$, and moreover, the diagram 
\begin{align*}
\xymatrix{ 
P^+ \ar^{\rho}[rrr]   & & &  T^+ \\
\Phi(\Omega^1_u(P)) \ar[u]^{\sigma} \ar_{\Phi(\mathrm{ver})}[rrr]& & &  \Phi(P \otimes T^+) \ar[u]_{\counit}
}
\end{align*}
commutes.
\end{enumerate} 
\end{lem} 
\begin{proof}
~~\\
\begin{enumerate}
\item Since $\mathrm{ver}$ is by construction a morphism in ${}^A_P\mathrm{Mod}$, we need only show that $\mathrm{ver}$ is a right $K$-comodule map, which is to say, we need to show that the following diagram is commutative:
        \begin{align*}
\xymatrix{ 
\Omega^1_u(P)   \ar[rrr]^{\mathrm{ver}}       \ar[d]_{\Delta_{R,K}}                  & & &    P \otimes T^+  \ar[d]^{\Delta_{R,K} \otimes \Ad_{R,K}} \\
\Omega^1_u(P) \otimes T   \,\,   \ar[rrr]_{\mathrm{ver} \otimes \mathrm{id}}  & & & P \otimes T^+ \otimes T.
}
\end{align*}
This follows from the calculation
    \begin{align*}
\xymatrix{ 
p\exd q    \ar[rrr]^{\mathrm{ver}}       \ar[d]_{\Delta_{R,K}}                  & & &    pq_{(1)} \otimes \rho(q^+_{(2)}) \ar[d]^{\Delta_{R,K} \otimes \Ad_{R,K}}   \\
p_{(1)}\exd q_{(1)} \otimes \rho(p_{(2)}q_{(2)})   \,\,   \ar[rrr]_{\mathrm{ver} \otimes \mathrm{id}}  & & & p_{(1)} q_{(1)} \otimes \rho(q^+_{(2)}) \otimes \rho(p_{(2)}q_{(3)}),
}
\end{align*}
where we have used the fact that 
\begin{align*}
(\Delta_{R,K} \otimes \Ad_{R,K})(pq_{(1)} \otimes \rho(q_{(2)})) = \, & p_{(1)}q_{(1)} \otimes \rho(q_{(4)}) \otimes \rho(p_{(2)}q_{(2)})\rho(S(q_{(3)})q_{(5)})\\
= &  \,  p_{(1)}q_{(1)} \otimes \rho(q_{(2)}) \otimes \rho(p_{(2)}q_{(3)}).
\end{align*}

\item The fact that $\sigma$ is an isomorphism in ${}^A_P\mathrm{Mod}$ is well known (see for example \cite[Theorem 2.10.3]{MMF2}.) Thus it remains to show that $\sigma$ is a right $K$-comodule map. This follows from the calculation 
    \begin{align*}
\xymatrix{ 
[\exd p]    \ar[rrr]^{\sigma}       \ar[d]_{\Delta_{R,K}}                  & & &    p^+ \ar[d]^{\Ad_{R,K}}   \\
[\exd p_{(2)}]\otimes \rho(S(p_{(1)})p_{(3)})   \,\,   \ar[rrr]_{\mathrm{\sigma} \otimes \mathrm{id}}  & & & [p^+_{(2)}]\otimes \rho(S(p_{(1)})p_{(3)}).
}
\end{align*}

\item Follows from
\begin{align*}
    \Psi(T^+) = A \square_{H} T^+ = A^{\co(H)} \otimes T^+ = P \otimes T^+.
\end{align*}
\item The calculation
    \begin{align*}
\xymatrix{ 
p^+   \ar[rrr]^{\rho}                     & & &    \rho(p)   \\
[\exd p]   \,\, \ar[u]^{\sigma}  \ar[rrr]_{\Phi(\mathrm{ver})}  & & & p_{(1)} \otimes \rho(p_{(2)}). \ar[u]_{\counit}
}
\end{align*}
now implies commutativity of the diagram as claimed. Moreover, since commutativity implies that $\rho|_{P^+}$ is a compositions of morphisms in  ${}^H\mathrm{Mod}^{K,\square}$, we can conclude that $\rho|_{P^+}$ is a morphism.
\end{enumerate}
\end{proof}

\begin{cor}
The diagram 
    \begin{align*}
\xymatrix{ 
\Omega^1_u(P)  \ar[d]^{(\mathrm{id} \otimes \sigma) \circ \unit} \ar[rrr]^{\mathrm{ver}}                     & & &  P \otimes T^+,   \\
A \square_H P^+  \,\,   \ar[rrru]_{\id \otimes \rho|_{P^+}}  & & &
}
\end{align*}  
commutes.
\end{cor}
\begin{proof} 
From \eqref{lem:PPver} we know the following diagram commutes
 \begin{align*}
\xymatrix{ 
\Omega^1_u(P)   \ar[rrr]^{\mathrm{ver}}       \ar[d]_{\mathrm{U}}                  & & &    P \otimes T^+ \ar[d]_{\mathrm{U}}  \\
A \square_H \Phi \left(\Omega^1_u(P)\right)  \,\, \ar[d]_{\mathrm{id}\otimes\sigma}\ar[rrr]_{\id \otimes \Phi(\mathrm{ver})}  & & & A \square_H \Phi( P \otimes T^+ )\ar[d]_{\mathrm{id}\otimes \mathrm{C}}
\\A \square_H P^+   \,\, \ar[rrr]_{\mathrm{id} \otimes \rho}  & & & P \otimes T^+  \ar@/^-4.0pc/[uu]_{\mathrm{id}_{P\otimes T^+}}.}
\end{align*}
which implies that the given diagram commutes. 
\end{proof}

With the above technical results about universal quantum principal bundles in hand, we now establish an existence result for non-universal principal bundles and principal pairs.

\begin{prop} \label{prop:PPQPB}
Let $\Omega^1(P)$ be a left $A$-covariant, right $K$-covariant, first-order differential calculus  over $P$, with corresponding sub-object  $N \sseq \Omega^1_u(P)$ in the category ${}^A_P\mathrm{Mod}^K$. It holds that 
\begin{enumerate}
\item the pair
$
(P,\Omega^1(P))
$
is a quantum principal bundle,
\item $\sigma\left(\Phi(\Omega^1(P)_{\mathrm{hor}})\right) = \iota'(V(B))P = \{[p] \,|\, p \in \ker(\rho|_{P^+})\}$.
\end{enumerate}
\end{prop}
\begin{proof}
~~~
\begin{enumerate}
\item Since $\Phi_P(P \otimes T^+)$ and $T^+$ are isomorphic as objects in ${}^H\mathrm{Mod}^{K,\square}$, we see that the right $H$-coaction on $\Phi_P(P \otimes T^+)$ is trivial. Thus since  $\Phi_P(\mathrm{ver}(N))$ is a subobject of $\Phi_P(P \otimes T^+)$, we have that 
$$
\mathrm{ver}(N) \simeq P \otimes \Phi_P(\mathrm{ver}(N)).
$$
In particular, we see that $\mathrm{ver}(N)$ is a free left $P$-module admitting a basis of left $A$-coinvariant elements. The only such subobject of $P \otimes T^+$ is of the form $P \otimes C$, for $C$ some right coideal of $T^+$.    
\item For any $p'(\exd b)p \in \Omega^1(P)_{\mathrm{hor}}$, observe that 
\begin{align*}
\sigma\!\left([p'\exd b p]\right) = \e(p)[b^+p] = \e(p)[b^+]p.
\end{align*}
This gives us the first claimed identity. 

Note next that the commutative diagram in Lemma \ref{lem:PPver} gives
\begin{align*}
\sigma\!\left(\Omega^1(P)_{\mathrm{hor}}\right) = & \, \sigma\!\left(\Phi_P(\ker(\mathrm{ver}))\right) = \ker(\rho), 
\end{align*}
establishing the second identity.
\end{enumerate}
\end{proof}

\subsection{Principal Pairs and Principal Connections}

In this subsection we investigate principal connections for quantum principal bundles coming from principal pairs. 

\begin{lem}
Let $(\pi':A \to H,\pi:A \to K)$ be a principal pair.
\begin{enumerate}
\item For $i:T^+ \to P^+$ a morphism splitting $\pi|_{P^+}:P^+ \to T^+$, a left $A$-covariant principal connection is given by 
$$
\Pi_i := \mathrm{U}^{-1}\Psi(\sigma^{-1} \circ i \circ \pi \circ \sigma)\mathrm{U}: \Omega^1_u(P) \to \Omega^1_u(P).
$$

\item This gives a bijective correspondence between left $A$-covariant principal connections $\Pi:\Omega^1_u(P) \to \Omega^1_u(P)$ and morphisms splitting $\pi|_{P^+}$.  
   
\item It holds that 
\begin{align*}
\sigma\!\left(\Phi(\Omega^1_u(P)_{\mathrm{ver}})\right) = \sigma\!\left(\Phi(\mathrm{im}(\Pi_i))\right) = i(T^+).
\end{align*}

\item The image, under $\sigma \circ \Phi$, of the decomposition of $\Omega^1_u(P)$ into horizontal and vertical forms, is given by the decomposition
\begin{align}\label{eqn:Pdecomposition}
P^+ = \ker(\pi) \oplus i(T^+)
\end{align}
in the category ${}^H\mathrm{Mod}^{K,\square}$.
\end{enumerate}
\end{lem}
\begin{proof}
~~~~
\begin{enumerate}
    \item Consider a left $A$-covariant principal connection $\Pi$. As recalled in \textsection \ref{subsection:QPB}, principal connections are equivalent to right $K$-comodule splittings 
$$
s:P\otimes T^+ \to \Omega^1_u(P)
$$
of the vertical projection $\mathrm{ver}$. Note that the connection associated to $s$ is left $A$-covariant if and only if $s$ is a left $A$-comodule map. Such maps correspond under $\Phi$ to splittings of 
$$
\Phi(\mathrm{ver}): \Phi(\Omega^1_u(P)) \to \Phi(P \otimes T^+).
$$
By Lemma \ref{lem:PPver}, such maps correspond to splittings of  
$
\pi|_{P^+}:P^+ \to T^+,
$
as claimed.
\item We have the commutative diagram
\begin{align*}
\xymatrix{ 
P^+ \ar@<-2pt>^{\substack{i\\~}}[rrrr] \ar@<-2pt>^{\;\sigma}[d]  & & & & \ar@<+2pt>_{\counit\;}[d] \ar@<-2pt>^{{\substack{~\\\pi}}}[llll] T^+ \\
\Phi(\Omega^1_u(P)) \ar@<-2pt>[u]^{\sigma^{-1}\;} \ar@<-2pt>_{\Phi(\mathrm{ver})}[rrrr]& & & & \Phi(P \otimes T^+) \ar@<+2pt>[u]_{\;\counit^{-1}} \ar@<-2pt>_{\Phi(s)}[llll]
}
\end{align*}
Thus we see that for $i$ a splitting of $\pi|_{P^+}$, the corresponding principal connection is explicitly given by $\Psi(\sigma^{-1} \circ i \circ \pi \circ \sigma)$. Moreover, part 1 above implies that 
\begin{align*}
  \sigma\!\left(\Phi(\Omega^1_u(P)_{\mathrm{ver}})\right) = & \, \sigma\!\left(\Phi(\mathrm{im}(\Pi_i))\right)
    =  \, \sigma\!\left(\mathrm{im}\!\left(\Phi(s \circ \mathrm{ver})\right)\right)
    =  \, \mathrm{im}\!\left(\sigma \circ \Phi(s) \circ \Phi(\mathrm{ver})\right).
\end{align*}
\item From the commutative diagram above, we see that 
$$
\sigma \circ \Phi(s) \circ \Phi(\mathrm{ver}) = \sigma^{-1} \circ \sigma \circ i \circ \counit \circ \counit^{-1} \circ \pi \circ \sigma = i \circ \pi \circ \sigma.
$$
Since $\pi$ and $\sigma$ are surjective maps, we have that 
$$
\mathrm{im}(i \circ \pi \circ \sigma) = \mathrm{im}(i) = i(T^+),
$$
giving the third claimed identity.
\item This follows directly from part {\em 3} of this proposition and from part {\em 2} of proposition \eqref{prop:PPQPB}.
\end{enumerate}
\end{proof}

Let us now consider a non-universal quantum principal bundle structure
$$
\Omega^1(P) = \Omega^1_u(P)/N
$$ 
for a principal pair $(\pi':A\to H,\pi:A \to K)$. As usual we denote $P:= A^{H}$ and $T = \pi(P)$.

\begin{prop}
A universal left $A$-covariant principal connection $\Pi$, corresponding to a splitting $i$ of the map $\pi|_{P^+}$, descends to a principal connection on $\Omega^1(P)$ if and only if 
$$
J := \sigma\left(\Phi(N)\right) 
$$ 
is a homogeneous subspace with respect to the decomposition in \eqref{eqn:Pdecomposition}.
\end{prop}
\begin{proof}
We note first that $\Pi_i$ descends to a connection on $\Omega^1(P)$ if and only if $N$ is homogeneous with respect to the decomposition of $\Omega^1_u(P)$ into horizontal and vertical forms. Since $\Phi$ is an equivalence it commutes with pullbacks, and so, it holds that 
\begin{align*}
\sigma\!\left(\Phi(N \cap P\Omega^1_u(B)P)\right) = \sigma\!\left(\Phi(N)\right) \cap \sigma\!\left(\Phi(P\Omega^1_u(B)P)\right),
\end{align*}
and moreover that
\begin{align*}
\sigma\!\left(\Phi(N \cap \mathrm{im}(\Pi_i)\right) = \sigma\left(\Phi(N)\right) \cap \sigma\!\left(\Phi(\mathrm{im}(\Pi_i))\right).  
\end{align*}
Thus we see  that $N$ is homogeneous with respect to the decomposition into horizontal and vertical forms if and only if $J$ is homogeneous with respect to the decomposition in \eqref{eqn:Pdecomposition}. Hence $\Pi_i$ will descend to a connection on $\Omega^1(P)$ if and only if $J$ is a homogeneous subspace.
\end{proof}

From this we get the following immediate corollary, which we present as such for sake of clarity.

\begin{cor}
If $\Pi$ descends to a principal connection on $\Omega^1(P)$, then 
\begin{align}
 V_{\mathrm{ver}}(P) :=   \sigma\!\left(\Phi(\Omega^1(P)_{\mathrm{ver}})\right) =  \{[i(t)] \,|\, t \in T^+\},
\end{align}
in the category ${}^H\mathrm{Mod}^{K,\square}$.
\end{cor}

\section{The Quantum Grassmannians $\OO_q(\mathrm{Gr}_{n,m})$} \label{section:quantumGrassPrelim}

In this section we introduce the main object of study in this paper, namely the quantum Grassmannians, their covariant line bundle, and the quantum Grassmannian sphere. In addition to defining these objects, we recall that a set of generators is given by the quantum minor determinants of $\OO_q(SU_n)$,  and recall how the coquasitriangular structure of $\OO_q(SU_n)$ acts on quantum minors. Since coquasi-triangular structures play a central role in the paper, we adopt an FRT approach throughout \cite[\textsection 9.1.1]{KSLeabh}.

\subsection{Quantum Matrices and Quantum Minors}

For $q \in \mathbb{R}_{>0}$, let $\OO_q(M_{n})$ be the FRT bialgebra associated to the $A$-series $R$-matrix
\begin{align}\label{eqn:Rmatrix}
R^{ij}_{kl} := q^{\delta_{ij}}\delta_{ik}\delta_{jl} + (q-q\inv)\theta(i-j)\delta_{il}\delta_{jk}, & & i,j,k,l = 1, \dots, n,
\end{align}
where $\theta$ is the Heaviside step function. As is well-known, the center of $\OO_q(M_{n})$ is generated by the grouplike element 
$$
\mathrm{det}_{n} := \sum\nolimits_{\sigma \in S_{n}}(-q)^{\ell(\sigma)}u^1_{\sigma(1)}u^2_{\sigma(2)} \cdots u^n_{\sigma(n)},
$$
which we call the  {\em quantum determinant}.

The quantum determinant is a special example of a quantum minor, whose general construction we now recall. Let $I: = \{i_1, \dots, i_m\}$ and $J:=\{j_1, \dots, j_p\}$ be a pair of subsets of $\{1,\dots, n\}$.
The associated {\em quantum  minor} $[I|J]$ is the element of $\OO_q(M_{n})$ given by
\begin{align*}
[I|J] := \sum_{\sigma \in S_p} (-q)^{\ell(\sigma)}u^{\sigma(i_1)}_{j_1} \cdots u^{\sigma(i_p)}_{j_p} =  \sum_{\sigma \in S_p} (-q)^{\ell(\sigma)} u^{i_1}_{\sigma(j_1)} \cdots u^{i_p}_{\sigma(j_p)}.
\end{align*} 
Note that when $I=J=\{1, \dots, n\}$ we get back the quantum determinant $\mathrm{det}_{n}$. The coproduct acts on any quantum minor $[I|J]$ according to 
\begin{align} 
\Delta([I|J]) = \sum_{K}  [I|K] \oby [K|J],
\end{align}
where summation is over all ordered subsets  $K \sseq \{1, \dots, n\}$ with $|K|=|I|$ \cite[\textsection 1.2]{NoumiUn}.

For $I,$ ordered subsets of $\{1,\dots ,n\}$ we define  $\ell(I,J):= |\{(i,j)\in I\times  \,| \, i>j\}|$ and denote by $I^{\perp}$ the ordered complement to $I$ in $\{1,\dots, n\}$. The following useful generalisation of Laplace expansions was established in \cite[Proposition 1.1]{NoumiUn}. For $I,J \sseq \{1, \dots, n\}$ with $|I| = |J|$, and $J_1$ a choice of non-empty subset of $J$, the associated \emph{Laplace expansion} is given by
\begin{align} \label{Laplace}
(-q)^{\ell(J_1,J_1^\perp)}  [I|J] = \sum_{I_1} (-q)^{\ell(I_1,I_1^\perp)}[I_1|J_1][I_1^\perp|J_1^\perp],
\end{align}
where summation is over all subsets $I_1 \sseq I$ such that  $|I_1| = |J_1|$. The $\ast$-map of $\OO_q(SU_n)$ acts on quantum minors as 
\begin{align} \label{eqn:starmin}
([I|J])^* = S([J|I]) = (-q)^{\ell(J,J^\perp) - \ell(I,I^\perp)}[{I^\perp}|{J^\perp}].
\end{align}

Recall that a {\em Young diagram}  is a finite collection of boxes arranged in left-justified rows, with the row lengths in non-increasing order. Young diagrams with $r$  rows are clearly equivalent to  {\em partitions} of {\em length} $r$, which is to say  elements 
\begin{align*}
\lambda = (\lambda_1, \ldots, \lambda_r) \in \mathbb{Z}^{r}_{>0}, & & \text{ such that ~} \lambda_1 \geq \cdots \geq \lambda_r.
\end{align*} 
We denote the set of partitions of length $r$ by $\mathcal{P}_r$. 
For $i = 1, \dots, r$, let  $e_i$ be the standard generators of the monoid  $\mathbb{Z}^r_{\geq 0}$. The {\em fundamental} partitions are given by
\begin{align*}
\varpi_k := e_1 + e_2 + \cdots + e_k, & & \text{ for } k= 1, \dots,  r.
\end{align*}
Note that any partition $\mu = (\mu_1,\dots, \mu_r)$ can be expressed as the sum of fundamental partitions: 
\begin{align*}
\mu = (\mu_1 - \mu_2) \varpi_1 + (\mu_2 - \mu_3) \varpi_2 + \cdots + (\mu_{r-1} - \mu_r)\varpi_{r-1} + \mu_r \varpi_r.
\end{align*}

A   {\em semi-standard tableau} of {\emph  shape} $\lambda$ with {\em labels} in $\{1,\dots, n\}$ is a collection  $\mathbf{T} = \{T_{a,b}\}_{a,b}$ of elements of $\{1,\dots, n\}$ indexed by the boxes of the corresponding Young diagram, and satisfying, whenever defined, the inequalities 
\begin{align*}
T_{a-1,b} < T_{a,b}, & & T_{a,b-1} \leq T_{a,b},
\end{align*} 
We denote by $\sstab(\lambda)$ the set of all semi-standard tableaux, for any dominant  weight $\lambda$.
\begin{eg} For the reader's convenience let us consider the simple example of the Young diagram corresponding to the partition $\{3,2,1\}$. The following collection of labels
\begin{center}
\ytableausetup{mathmode, boxsize=1.5em}\begin{ytableau}
a_1 & b_1 & c_1 
\\a_2 & b_2
\\a_3 \end{ytableau}
\end{center}
is a Young tableau if and only if $a_1<a_2<a<3$, $b_1<b_2$, and  $a_1\leq b_1\leq c_1$, $a_1\leq b_1$.
\end{eg}  

For a Young diagram with $n-1$ rows, and a semi-standard tableau $\mathbf{T}$, we associate the product 
\begin{align*}
[{\mathbf{T}}] := [T_1|(T_1)] \cdots [T_{\lambda_1}|(T_{\lambda_1})] \in \OO_q(M_{n}),
\end{align*}
where $T_b := \{T_{1,b}, \ldots, T_{l_b,b}\}$, and $(T_b) := \{1, \dots, b\}$, as ordered sets, for $1 \leq b \leq \lambda_1$, and $l_b$ is the length of the $b^{\mathrm{th}}$ column. Associated to any $\lambda$ we have the coideal
\begin{align*}
V(\lambda) := \mathrm{span}_{\mathbb{C}} \{[\mathbf{T}]  \, | \, \mathbf{T} \in \sstab(\lambda)\}.
\end{align*}
We denote the associated subcoalgebra of matrix coefficients by $C(\lambda)$. (Note that in the definition of $V(\lambda)$ taking $(T_b) = \{1, \dots, b\}$ is merely a convention. One could just as well choose $(T_b) = \{n-b+1, \dots, n\}$, and produce an isomorphic $V(\lambda)$. By abuse of notation we use the same symbol for both comodules, the sense always being clear from the context.)

\subsection{The Hopf Algebras $\OO_q(SU_n)$ and $\OO_q(U_n)$}

Centrality of $\mathrm{det}_{n}$ makes it easy to adjoin a formal inverse $\det_{n}^{-1}$ to $\OO_q(M_{n})$. The bialgebra structure of $\OO_q(M_{n})$ uniquely extends to a bialgebra structure on this extended algebra. Indeed it is a Hopf algebra with antipode defined by 
\begin{align*}
S(\mathrm{det}_{n} \inv) := \dt_{n}, ~~~~ S(u^i_j) := (-q)^{i-j}\sum\nolimits_{\sigma \in S_{n-1}}(-q)^{\ell(\sigma)}u^{k_1}_{\sigma(l_1)}u^{k_2}_{\sigma(l_2)} \cdots u^{k_{n-1}}_{\sigma(l_{n-1})}\dt_n\inv,
\end{align*}
where $\{k_1, \ldots ,k_{n-1}\} := \{1, \ldots, n\}\bs \{j\}$, and $\{l_1, \ldots ,l_{n-1}\} := \{1, \ldots, n\}\bs \{i\}$ as ordered sets. A Hopf $*$-algebra structure is determined by ${(\mathrm{det}_{n}^{-1})^* = \mathrm{det}_{n}}$, and $(u^i_j)^* =  S(u^j_i)$. We denote this Hopf $*$-algebra by $\OO_q(U_{n})$, and call it the {\em quantum unitary group of order $n$}. We denote the Hopf  {$*$-algebra} $\OO_q(U_{n})/\langle \mathrm{det}_{n} - 1 \rangle$ by $\OO_q(SU_{n})$, and call it the {\em quantum special unitary group of order $n$}. The $\ast$-map of $\OO_q(SU_n)$ acts on quantum minors as 
\begin{align} \label{eqn:starmin}
([I|J])^* = S([J|I]) = (-q)^{\ell(J,J^\perp) - \ell(I,I^\perp)}[{I^\perp}|{J^\perp}].
\end{align}
An analogous formula holds for $\OO_q(U_n)$.

The coalgebras $C(\lambda)$ embed into $\OO_q(U_n)$ and $\OO_q(SU_n)$, and by abuse of notation we denote them by the same symbol. Both Hopf algebras $\OO_q(SU_{n})$ and $\OO_q(U_{n})$ are \emph{cosemisimple}, which is to say, they are the direct sum of their simple subcoalgebras. Explicitly, we have
\begin{align} \label{eqn:PW}
    \OO_q(U_{n}) \simeq \bigoplus_{\lambda \in\mathcal{P}_{n-1}} \bigoplus_{k \in \mathbb{Z}} C(\lambda)\mathrm{det}^k_{n}, & &  \OO_q(SU_{n}) \simeq \bigoplus_{\lambda \in \mathcal{P}_{n-1}} C(\lambda),
\end{align}
as established in \cite[Theorem 2.11]{noumi0} and \cite{NoumiUn}. We denote the associated Haar maps $\haar_{\OO_q(U_n)}$ and $\haar_{\OO_q(SU_n)}$. (Recall that $\haar_{\OO_q(U_n)}$ and $\haar_{\OO_q(SU_n)}$ are given by projection onto the subcoalgebra $\mathbb{C}1_{\OO_q(U_n)}$, and $\mathbb{C}1_{\OO_q(U_n)}$ respectively.) Both maps are positive, that is   
\begin{align*}
\haar_{\OO_q(U_n)}(a^*a) > 0, & & \text{ for all non-zero } a \in \OO_q(U_n), 
\end{align*}
and similarily for $\haar_{\OO_q(SU_n)}$. In particular, this means that each Hopf algebra can be completed to a compact quantum group in the sense of Woronowicz \cite{KoornDijk}.

\subsection{Quantum Grassmannians}

We begin by recalling Meyer's definition of the quantum Grassmannians \cite{MEYER}. (As discussed in \textsection \ref{Appendix:GensofGrassSphere} this is equivalent to Stokmann and Dijkhuizen's definition given \cite{DijkStok}.)  Let $\alpha_m:\OO_q(SU_{n}) \to \OO_q(U_m)$ be the surjective Hopf $*$-algebra map defined on generators by 
\begin{align*}
 &\alpha_m(u^i_j) := \delta_{ij}\dt^{-\delta_{in}}_{m},&& \text{ for } i,j = 1, \dots, n; ~ (i,j) \notin  M \by M,\\
 & \alpha_m(u^i_j) :=  u^{i}_{j},&& \text{ for } (i,j) \in M \by M. ~~~~~~~~~~~~~~~~~~~~~
\end{align*}
Moreover, let $\beta_m:\OO_q(SU_{n}) \to  \OO_q(SU_{n-m})$ be the following surjective Hopf \\$*$-algebra map 
\begin{align*}
& \beta_m(u^i_j) := \delta_{ij} 1, \,\,\quad\qquad && \text{ for } i,j = 1, \dots, n; ~ (i,j) \notin  M^\perp \by M^\perp,\\
&  \beta_m(u^i_j) :=  u^{i-m}_{j-m}, \,\,\quad\qquad & &  \text{ for } (i,j) \in M^\perp \by M^\perp.
\end{align*}

\begin{defn}
Denoting $\OO_q(L_m) := \OO_q(U_m) \oby \OO_q(SU_{n-m})$, the {\em quantum Grassmannian} $\OO_q(\mathrm{Gr}_{n,m})$ is the quantum homogeneous space associated to the surjective Hopf $*$-algebra surjection 
\begin{align*}
\pi:\OO_q(SU_{n}) \to \OO_q(L_m), & & \pi:= (\alpha_m \oby \beta_m) \circ \Delta.
\end{align*}
\end{defn}
Since $\OO_q(L_m)$ is the product of two cosemisimple Hopf algebras, it is itself cosemisimple. Hence, as discussed in \textsection \ref{subsection:PCA}, the pair $(\OO_q(SU_n),\Delta_{R,\OO_q(L_m)})$ is a principal comodule algebra. 

Classically the Grassmannians $\mathrm{Gr}_{n,m}$ are projective manifolds, with an explicit embedding of $\mathrm{Gr}_{n,m}$ into complex projective space being given by the venerable Pl\"ucker embedding. Just as for the coordinate algebra $\OO(\mathrm{Gr}_{n,m})$, the  homogeneous coordinate ring $S(\mathrm{Gr}_{n,m})$ associated to the Pl\"ucker embedding admits a natural $q$-deformation.

\begin{defn}
The \emph{twisted Grassmannian homogeneous coordinate ring}  $S_q(\text{Gr}_{n,m})$  is  the subalgebra of $\OO_q(SU_{n})$ generated by the \emph{quantum minors} $\overline{z}_I$.
\end{defn}

As suggested by the notation, the algebra $S_q(\mathrm{Gr}_{n,m})$ does indeed reduce to the homogeneous coordinate ring $S(\mathrm{Gr}_{n,m})$ when $q=1$. The ring theoretic properties of the twisted Grassmannian homogeneous coordinate rings have been extensively studied in the literature \cite{Fioresi2004, KLR2004,LR2006,LLR2008}. They form a very important family of examples in the theory of quantum cluster algebras \cite{GrL2011,GrL2014}. Moreover, as Noetherian connected $\mathbb{Z}_{\geq 0}$-graded algebras, with well-behaved homological regularity, they are important examples of noncommutative projective varieties. In fact, the special case of quantum projective space $S_q(\mathbb{CP}^{n})$, where the homogeneous coordinate ring reduces to the well-known quantum plane, is a prototypical example in noncommutative projective geometry \cite{ATvdB,Avdb}.

\begin{remark}
In the literature on twisted Grassmannian coordinate algebras, $S_q(\mathrm{Gr}_{n,m})$ is usually defined to be the subalgebra of the quantum matrices $\OO_q(M_n)$ generated by the quantum minors, where quantum minors in $\OO_q(M_n)$ are defined in obvious analogy with the $\OO_q(SU_n)$ case. However, the canonical projection from $\OO_q(M_n)$ to $\OO_q(SU_n)$ restricts to an algebra embedding of $S_q(\mathrm{Gr}_{n,m})$ allowing us to identify the two presentations.
\end{remark}

\begin{remark}
Classically the Grassmannians are very special examples of generalised flag manifolds $G/L_S$, where $G$ is a compact semisimple Lie group, and $L$ is a Levi subgroup of $G$. This picture generalises to the quantum setting. Indeed for every generalised flag manifold $G/L_S$ we have $q$-deformed coordinate algebra $\OO_q(G/L_S)$ \cite{DijkStok} and a $q$-deformed homogeneous coordinate ring \cite{LR91, Soib}. Indeed the present paper should be considered as part of a general program to understand the noncommutative complex, K\"ahler, and projective geometry of the quantum flag manifolds.
\end{remark}

\subsection{Line Bundles and the Quantum Grassmann Sphere}

We now introduce a variation on the map $\alpha_m$. Let $\alpha'_m:\OO_q(SU_{n}) \to \OO_q(SU_m)$ be the surjective Hopf $*$-algebra map defined on generators by 
\begin{align*}
 &\alpha'_m(u^i_j) := \delta_{ij}1,~~~~~~~~~~~~~~~~~~\,~~ \text{ for } i,j = 1, \dots, n; ~ (i,j) \notin  M \by M,\\
  &\alpha'_m(u^i_j) :=  u^{i}_{j},~~~~~ ~~~\,~~~~~~~~~~~~~~ \text{ for } (i,j) \in M \by M. ~~~~~~~~~~~~~~~~~~~~~
\end{align*}
Associated to $\alpha'_m$ we have a quantum homogeneous space, which as we see below, gives us a very useful description of the covariant line bundles over the quantum Grassmannians.

\begin{defn}
Denoting $\OO_q( L^{\mathrm{s}}_m) := \OO_q(SU_m) \oby \OO_q(SU_{n-m})$, the {\em quantum Grassmann sphere} $\OO_q(S^{n,m})$ is the quantum homogeneous space associated to the surjective Hopf $*$-algebra surjection 
\begin{align*}
\pi':\OO_q(SU_{n}) \to \OO_q(L^{\mathrm{s}}_m), & & \pi':= (\alpha'_m \oby \beta_m) \circ \Delta.
\end{align*}
\end{defn}

Just as for the quantum Grassmannians, the pair $(\OO_q(SU_n),\Delta_{R,\OO_q(L^{\mathrm{s}}_m)})$ is a principal comodule algebra. 

\begin{lem}
It holds that 
\begin{enumerate}
\item $\OO_q(S^{n,m})$ is a right $\OO_q(L_m)$-subcomodule of $\OO_q(SU_n)$,
\item $\pi(\OO_q(S^{n,m})) = \langle \mathrm{det}_{m}^{\pm 1} \otimes 1 \rangle \simeq \OO(U_1)$.
\end{enumerate}
\end{lem}
\begin{proof}
Take the decomposition of $\OO_q(SU_n)$ into irreducible right $\OO_q(L_m)$-subcomodules:
\begin{align} \label{eqn:leftKDecomp}
\OO_q(SU_n) \simeq \bigoplus_{\mu \in Y} V_{\mu},
\end{align}
where $Y$ is some subset of index set of isomorphism classes of $\OO_q(L_m)$-comodules. Now for any $\mu$, we denote by $C(V_{\mu})$ the coefficient coalgebra of $V_{\mu}$, which is necessarily simple. 
Now the Peter--Weyl decompositions given in \eqref{eqn:PW} imply that the canonical projection $\mathrm{proj}:\OO_q(U_m) \to \OO_q(SU_m)$ is injective when restricted to a simple subcoalgebra, and moreover that the image under $\mathrm{proj}$ of a simple subcoalgebra is again simple. Since  
$$
\Delta_{R,\OO_q(L^{\mathrm{s}}_m)}(V_{\mu}) \sseq V_{\mu} \otimes \mathrm{proj}(C(V_{\mu})),
$$
we now see that $V_{\mu}$ must be irreducible as a right $\OO_q(L^{\mathrm{s}}_m)$-comodule. In particular, $V_{\mu}$ can contain a right $\OO_q(L^{\mathrm{s}}_m)$-coinvariant element only if $\mathrm{dim}(C(V_{\mu})) = 1$, and hence only if 
\begin{align*}
C(V_{\mu}) = \mathrm{det}_m^k \otimes 1, & & \textrm{ for some } k \in \mathbb{Z}.
\end{align*}
Thus since $\OO_q(S^{n,m})$ is clearly homogeneous with respect to the decomposition in \eqref{eqn:leftKDecomp}, we must have that
$$
\Delta_{R,\OO_q(L^{\mathrm{s}}_m)}(\OO_q(S^{n,m})) \sseq \OO_q(S^{n,m}) \otimes \langle \mathrm{det}^{\pm 1}_m \otimes  1 \rangle.
$$
It now follows from Lemma \ref{lem:PPsetup} that 
$$
\pi'(S^{n,m}) \sseq \langle \mathrm{det}^{\pm 1}_m  \otimes 1\rangle.
$$
Finally, the definition of $\pi$ implies that 
\begin{align*}
\pi(z) = \mathrm{det}_m \otimes 1, & & \pi(\overline{z}) = \mathrm{det}^{-1}_m \otimes 1,
\end{align*}
giving us the opposite inclusion and hence equality.
\end{proof}

As explained in \cite[Example 1.6.7]{Monty}, the fact that $\OO(U_1)$ is the group Hopf algebra of $\mathbb{Z}$, means that an associated $\mathbb{Z}$-graded algebra structure
$$
\OO_q(S^{n,m}) \simeq \bigoplus_{k \in \mathbb{Z}} \EE_k.
$$
The following proposition presents some properties of this grading which will be used throughout the paper.

\begin{prop}
Denoting by $\Delta_{R,\OO(U_1)}$ the right $\OO(U_1)$-coaction corresponding to the $\mathbb{Z}$-algebra grading of $\OO_q(S^{n,m})$, the pair $(\OO_q(S^{n,m}),\Delta_{R,\OO(U_1)})$ is a $\OO(U_1)$-principal comodule algebra. Moreover, each $\EE_k$ is a covariant line bundle over $\EE_0 = \OO_q(\mathrm{Gr}_{n,m})$.
\end{prop}
\begin{proof}
By definition 
$$
\EE_0 = \OO_q(S^{n,m})^{\co(\OO_q(L_m))} = \OO_q(SU_n)^{\co(\OO_q(L_m))} = \OO_q(\mathrm{Gr}_{n,m}),
$$
where the second identity follows from the fact that all right $\OO_q(L_m)$-coinvariant elements of $\OO_q(SU_n)$ are automatically contained in $\OO_q(S^{n,m})$. By construction each $\EE_k$ is a left $\OO_q(SU_n)$-comodule, and so it is a relative Hopf module over $\OO_q(\mathrm{Gr}_{n,m})$. 

The fact that we have an algebra grading means that $\EE_1\EE_{-1} \sseq \EE_{0} = \OO_q(\mathrm{Gr}_{n,m})$, which is to say,  $\EE_1\EE_{-1}$ is a subobject of an irreducible object. Since there are no zero divisors in $\OO_q(SU_n)$, it cannot be that $\EE_1\EE_{-1} = 0$, hence we must have $\EE_1\EE_{-1} = \EE_0$ implying that the grading is strong, and hence that we have a Hopf--Galois extension (see for example \cite[Theorem 8.1.7]{Monty}).

The fact that the grading is strong in turn implies that 
\begin{align*}
\EE_k \otimes_{\OO_q(\mathrm{Gr}_{n,m})} \EE_l = \EE_k \otimes_{\EE_0} \EE_l \simeq \EE_k\EE_l, & & \textrm{ for all ~} k,l \in \mathbb{Z},
\end{align*}
see \cite[Corollary 3.1.2]{FreddyO}. In particular, for all $k \in \mathbb{Z}$, we have
$$
\EE_k\otimes_{\OO_q(\mathrm{Gr}_{n,m})} \EE_{-k} \simeq  \EE_{-k} \otimes_{\OO_q(\mathrm{Gr}_{n,m})} \EE_{k} \simeq \OO_q(\mathrm{Gr}_{n,m}),
$$
which is to say, each $\EE_k$ is a line bundle. 

Note next that since $\EE_k$ is invertible, it is projective as a left and as a right module over $\OO_q(\mathrm{Gr}_{n,m})$. This means that $\OO_q(S^{n,m})$ is a direct sum of projective modules and hence also projective as a left and as a right module, and hence flat as a left and as a right module. Finally, we note that since $\OO_q(S^{n,m})$ contains $\OO_q(\mathrm{Gr}_{n,m})$ as a direct summand, it must be faithfully flat as a right and left module. Thus we have a principal comodule algebra as required.
\end{proof}

Finally, collecting all the results and observations of this subsection allows to conclude the following corollary.

\begin{cor}
The pair of Hopf algebra maps $(\pi,\pi')$ is a principal pair.
\end{cor}

\subsection{Sets of Generators}

As discussed in the appendix, the algebra $\OO_q(\mathrm{Gr}_{n,m})$ is generated by the quantum minor products
\begin{align*}
    z^{IJ} := z_I\overline{z}_J & & \textrm{ for } |I| = m, \, |J| = n-m.
\end{align*}
Using these generators it is possible to produce a set of generators for $\OO_q(S^{n,m})$, as well as a very useful presentation of the line bundles introduced in the previous subsection. We begin with an abstract lemma used in the proof of the main proposition.

\begin{lem} \label{lem:TakGens}
Let $A$,and $H$ be two Hopf algebras, $\pi:A \to H$ a Hopf algebra map, and $B = A^{\co(H)}$ the associated quantum homogeneous space. Any irreducible $\F \in {}^A_B\mathrm{Mod}$ is generated as a left $B$-module by any left $A$-subcomodule of $\F$. 
\end{lem}
\begin{proof}
For $V$ a non-zero left $A$-subcomodule of $\F$, consider the left $B$-submodule $BV$ that it generates. Since $BV$ is clearly also a left $A$-subcomodule of $\F$, it is a sub-object of $\F$. However, since $\F$ is by assumption irreducible, and $BV \neq 0$, we must have that $BV = \F$.
\end{proof}

\begin{prop}
It holds that:
\begin{enumerate}
\item the quantum Grassmann sphere $\OO_q(S^{n,m})$ is generated as an algebra by the elements 
\begin{align} \label{eqn:GsphereGens}
\Big\{ z_I, \, \overline{z}_J  \,|\, \textrm{ for } |I| = n, |J| = n-m \Big\},
\end{align}
\item $V_{k\varpi_m}$ is a subcomodule of $\EE_k$, and  so, it generates $\EE_k$ as a left $\OO_q(\mathrm{Gr}_{n,m})$-module,
\item the $\mathbb{Z}$-algebra grading on $\OO_q(S^{n,m})$ is determined by 
\begin{align} \label{eqn:degree}
\mathrm{deg}(z_I) = - 1, & & \mathrm{deg}(\overline{z}_J) = 1.
\end{align}
\end{enumerate}
\end{prop}
\begin{proof}
 It follows from the definition of $\pi$ that 
\begin{align*}
    \Delta_R(z_I) = & \, \sum_{A,B}[I|A]\oby \alpha([A|B])\oby \beta(z_B)
     =  \, z_I \oby \alpha(z)\oby \beta(z)
  =  \, z_I\otimes \mathrm{det}_m \otimes 1,
\end{align*}
where the last equality follows from the fact 
$
\alpha(z) = \mathrm{det}_m \in \OO_q(U_m).
$
A similar computation shows that
\begin{align*} 
    \Delta_R(\overline{z}_J) = \overline{z}_I\otimes \mathrm{det}^{-1}_m \otimes 1.
\end{align*}
Thus the algebra generated by the elements given in \eqref{eqn:GsphereGens} is contained in $\OO_q(S^{n,m})$, and that $\mathrm{deg}(z_I) = -1$ and $\mathrm{deg}(\overline{z}_J) = +1$. 

As a direct implication we see that the irreducible comodule $V_{k\varpi_m}$ is contained in $\EE_k$. Hence Lemma \ref{lem:TakGens} implies that $V_{k\varpi_m}$ generates $\EE_k$ a left $\OO_q(\mathrm{Gr}_{n,m})$-module. Moreover, since $\OO_q(S^{n,m})$ is equal to the direct sum of all line bundles $\EE_k$, we see that $\OO_q(S^{n,m})$ is generated as a $\OO_q(\mathrm{Gr}_{n,m})$-module by elements given in $2$. But since $\OO_q(\mathrm{Gr}_{n,m})$ is generated by the elements $z^{IJ} = z_I\overline{z}_J$, we see that $\OO_q(S^{n,m})$ is in fact generated by the elements given in $\eqref{eqn:GsphereGens}$. Moreover, we see that knowing the degrees of $z_I$ and $\overline{z}_J$ completely determines the $\mathbb{Z}$-grading of $\OO_q(S^{n,m})$.
\end{proof}

We finish by giving an explicit construction of an $\ell$-map for the principal comodule algebra $(\OO_q(S^{n,m}),\Delta_{R,\OO(U_1)})$. This is a direct  direction of the construction of the $\ell$-map for the quantum Hopf fibration \cite[Example 24.4]{MajLeabh}.

\begin{prop}\label{prop:SCB}
A principal $\ell$-map $\ell: \mathcal{O}(U_1) \to \mathcal{O}_q(S^{n,m}) \oby \mathcal{O}_q(S^{n,m})$  is defined by 
\begin{align*}
t^k \mapsto S\big(z^k_{(1)}\big) \oby z^k_{(2)}, & & t^{-k} \mapsto S\big(\overline{z}_{(1)}^k\big) \oby \overline{z}_{(2)}^k, & & k \in \bZ_{\geq 0}.
\end{align*}
Hence,  $\mathcal{O}_q(\mathrm{Gr}_{n,m}) \hookrightarrow \mathcal{O}_q(S^{n,m})$ is a principal $\mathcal{O}(U_1)$-comodule algebra, which we call the {\em standard circle bundle} of $\mathcal{O}_q(\mathrm{Gr}_{n,m})$.
\end{prop}
\begin{proof}
We begin by showing that $\ell$ is well defined, which is to say that  its image is contained in $\mathcal{O}_q(S^{n,m}) \oby \mathcal{O}_q(S^{n,m})$. Note first that
\begin{align} \label{ellcoprod}
S\left(z_{(1)}\right) \oby z_{(2)} =   \sum_{K_1, \dots, K_k} S\left([M|K_k]\right) \cdots S\left([M|K_1]\right) \oby z_{K_1} \cdots z_{K_k}.
\end{align}
It follows from \eqref{eqn:starmin} that $\ell(t^k)$ is contained  in $\mathcal{O}_q(S^{n,m}) \oby \mathcal{O}_q(S^{n,m})$, for  $k \geq 0$, as required. The case for $k < 0$ follows similarly.
\\We now show that the requirements of definition \eqref{defn:ell} are satisfied. It is obvious that  conditions {\em 1} and {\em 2} hold. For $k\geq 0$,  condition {\em 3} follows from  
\begin{align*}
    (\ell \oby \id_{\mathcal{O}(U_1)}) \circ \DEL_{\mathcal{O}(U_1)}(t^k) = & \, \ell(t^k) \oby t^k \\
=  & S\left(z_{(1)}\right)^k \oby z_{(2)}^k \oby t^k \\
= & \, (\id_{\mathcal{O}_q[SU_n]} \oby \DEL_R)S(z_{(1)} \oby z_{(2)}\big)\\
= & \, (\id_{\mathcal{O}_q[SU_n]} \oby \DEL_R)\big(\ell(t^k)\big).
\end{align*}
The fourth condition is demonstrated analogously, as are both conditions for the case of $k<0$. 
\end{proof}

 \subsection{A Coquasitriangular Structure and the Goodearl Formulae}
 
 We begin by recalling the FRT-construction of a coquasitriangular structure for $\OO_q(SU_{n})$. (For details on general FRT-algebras see  \cite[\textsection 10.1.2]{KSLeabh}, and for the special case of $\OO_q(SU_n)$ see \cite[9.2]{KSLeabh}.) The $R$-matrix presented in \eqref{eqn:Rmatrix} is a solution of the quantum Yang--Baxter equation, thus by the general theory of FRT-algebras a coquasi-triangular structure for $\OO_q(SU_n)$ is defined by
\begin{align*}
r(u^i_j \otimes u^k_l) := \lambda R^{ik}_{jl},
\end{align*}
where $\lambda \in \mathbb{R}_{>0}$ is positive $n^{\mathrm{th}}$-root of $q^{-1}$.

We next recall a number of formulae, due to Goodearl \cite{Goodearl2005}, which describe the action of the coquasi-triangular structure of $\OO_q(SU_{n})$ on the quantum minors. We find it convenient to use the following notation: For $I \sseq \{1, \dots, n\}$, and $i,j \in \{1, \dots, n\}$,
\begin{align*}
I_{ij} :=
\left\{
	\begin{array}{ll}
		\big(I\bs\{i\}\big) \cup \{j\}  & \mbox{if } i \neq j \mbox{ and } i \in I \\
		I & \mbox{if }  i = j.
		 	\end{array}
\right.
\end{align*}
For $i \neq j$,  when $i \notin I$ or $j \in I$ we do not assign a direct meaning to the symbol $I_{ij}$, however we do denote  $[I_{ij}|J] = [J|I_{ij}] := 0,$ for all index sets $J$. The following result follows directly from the formulae established for the quantum matrices $\OO_q(M_{n})$ in \cite[\textsection 2]{Goodearl2005}.

\begin{prop}[Goodearl formulae]  It holds that, for $i, j = 1,\dots, n$,
\begin{enumerate}

\item $r(u^i_j\otimes [I|J]) \neq 0$ if and only if  $i \geq j$ \text{ and } $J = I_{ji}$,

\item $r([I|J]\otimes u^i_j) \neq 0$ if and only if $i \leq j$ and $J = I_{ji}$.

\end{enumerate}
\end{prop}
We inform the interested reader that a comprehensive  presentation of the action of the coquasitriangular structure of the quantum matrices $\OO_q(M_{n})$ on their quantum minors can be found in \cite{Goodearl2005}.

\subsection{The Heckenberger--Kolb Calculi and Holomorphic Structures}

In this subsection we give a brief presentation of the Heckenberger--Kolb calculi of the quantum Grassmannians. One should note that each of the results presented below extends directly to the more general setting of the irreducible quantum flag manifolds. All references refer to this more general form. 

We start  with the classification of first-order differential calculi over $\OO_q(\mathrm{Gr}_{n,m})$. Recall that a first-order calculus over an algebra $B$ is called \emph{irreducible} if it contains no non-trivial $B$-sub-bimodules. 

\begin{thm} \cite[\textsection 2]{HK} \label{thm:HKClass}
There exist exactly two, non-isomorphic, finite-dimensional irreducible left-covariant first-order differential calculi over $\OO_q(\mathrm{Gr}_{n,m})$. 
\end{thm}

We denote the direct sum of these two first-order calculi by $\Omega^1_q(\mathrm{Gr}_{n,m})$. Moreover, we denote the maximal prolongation of $\Omega^1_q(\mathrm{Gr}_{n,m})$ by 
$$
\Omega^\bullet_q(\mathrm{Gr}_{n,m}) = \bigoplus_{k \in \mathbb{Z}_{\geq 0}} \Omega^k_q(\mathrm{Gr}_{n,m}),
$$
and call it the \emph{Heckenberger--Kolb} calculus of the quantum Grassmannians.

\begin{thm}\label{CalcBasis} \cite[Proposition 3.11]{HKdR} Each $\Omega^k_q(\mathrm{Gr}_{n,m})$ has classical dimension, which is to say, 
\begin{align*}
\dim\!\big(\Phi(\Om^k_q(\mathrm{Gr}_{n,m}))\big) = \binom{2m(n-m)}{k}, & & \text{~ for ~}k = 0, \dots, 2m(n-m),
\end{align*}
and $\Om^k_q(\mathrm{Gr}_{n,m}) = 0$, for $k > 2m(n-m)$. 
\end{thm}

The calculus $\Omega^{\bullet}_q(\mathrm{Gr}_{n,m})$  has many remarkable properties, some of which we will now recall. Firstly, the results of~\cite{HK}, \cite{HKdR}, and~\cite{MarcoConj} give the following result about covariant complex structures.  

\begin{prop} \label{prop:complexstructure}
It holds that:
\begin{enumerate}
\item $\Omega^{\bullet}_q(\mathrm{Gr}_{n,m})$  is a $*$-calculus,
\item $\Omega^{\bullet}_q(\mathrm{Gr}_{n,m})$ admits precisely two left $\OO_q(SU_n)$-covariant complex structures, and these are opposite to the other. 
\end{enumerate}
\end{prop}

In the classical setting Liouville's theorem for a compact complex manifold $X$ says that the only holomorphic functions from $X$ to $\mathbb{C}$ are constant. The following result shows that this fact extends to the quantum setting. (A proof for the special case of quantum projective space can be found in \cite{KLVSCP1, KKCP2, KKCPN}, while the general irreducible quantum flag manifold case is treated in the forthcoming paper \cite{DOW}.)

\begin{thm}[Liouville]\label{theorem:Liouville}
The Heckenberger--Kolb calculus $\Omega^1_q(\mathrm{Gr}_{n,m})$ is connected, which is to say 
$$
H^0_{\adel,q}(\mathrm{Gr}_{n,m}) = H^0_{\del,q}(\mathrm{Gr}_{n,m}) = \mathbb{C}1.
$$
\end{thm}

Next we recall an existence and uniqueness result for covariant connections and holomorphic structures for relative Hopf modules over the the quantum Grassmannians. 
 
\begin{thm}\cite[Theorem 4.5]{HVBQFM} \label{thm:covconnection}
For each $\F \in {}^{~~\OO_q(SU_{n})}_{\OO_q(\mathrm{Gr}_{n,m})}\mathrm{Mod}_0$, it holds that 
\begin{enumerate}
\item $\F$ admits a left $\OO_q(SU_n)$-covariant connection 
$$
\nabla: \F \to \Omega^1_q(\mathrm{Gr}_{n,m}) \otimes_{\OO_q(\mathrm{Gr}_{n,m})} \F,
$$ 
and this is the unique such connection if $\F$ is simple,

\item a left $\OO_q(SU_n)$-covariant holomorphic structure for $\F$ is given by 
$$
\adel_{\F} := \mathrm{proj}^{(0,1)} \circ \nabla
$$
and this is the unique such holomorphic structure if $\F$ is simple.
\end{enumerate}
\end{thm}

We finish with the novel observation that the $*$-calculus structure of $\Omega^1_q(\mathrm{Gr}_{n,m})$ can be deduced from the classification of Heckenberger and Kolb without any need for direct calculation. As is easily seen, the argument is directly extendable to the general irreducible quantum flag manifolds case, even if we do not state it in this form.

\begin{prop}
The Heckenberger--Kolb calculus $\Omega^{\bullet}_q(\mathrm{Gr}_{n,m})$ is a differential $*$-calculus.
\end{prop}
\begin{proof}
Denote by $N^{(1,0)}$ the left $\OO_q(SU_n)$-covariant sub-bimodule of the universal calculus corresponding to $\Omega^{(1,0)}$. Since the universal calculus of any $*$-algebra is a first-order differential $*$-calculus, we can consider 
$$
W := (N^{(1,0)})^*.
$$
As is easily checked, $W$ is a left $\OO_q(SU_n)$-covariant sub-bimodule of the universal calculus. Moreover, since $\Omega^{(1,0)}$ is irreducible, the associated first-order calculus
$$
\Gamma^1 := \Omega^1_{u,q}(\mathrm{Gr}_{n,m}) 
$$
must also be irreducible. It now follows from the classification of irreducible  finite-dimensional left-covariant  calculi over the quantum Grassmannians \ref{thm:HKClass} that $\Gamma^1$ is equal to $\Omega^{(1,0)}$ or $\Omega^{(0,1)}$. In the first case $\Omega^{(1,0)}$ is necessarily a $*$-calculus, for the second case the Heckenberger--Kolb calculus $\Omega^1_q(\mathrm{Gr}_{n,m})$ is a $*$-calculus. Repeating the argument for $\Omega^{(1,0)}$ now implies that $\Omega^1_q(\mathrm{Gr}_{n,m})$ is a $*$-calculus as claimed.
\end{proof}

\begin{remark}
Each of the relative Hopf modules $\F$ comes endowed with a covariant  $q$-deformed Hermitian metric $g:\F \times \F \to \OO_q(\mathrm{Gr}_{n,m})$, which is unique when $\F$ is simple. Using $g$ one can mimic the classical construction of the Chern connection \cite[\textsection 8.6]{BeggsMajid:Leabh}, which by uniqueness must coincide with $\nabla$. Chern connections will be appear in \ref{section:Kahler} when we discuss positivity for line bundles. 
\end{remark}


\section{A Restriction Calculus Presentation of the Heckenberger--Kolb Calculus} \label{section:RestPresHK}

In this section we generalise the work of \cite{MMF1} for the case of quantum projective space and realise the Heckenberger--Kolb calculus of $\OO_q(\mathrm{Gr}_{n,m})$ as the restriction of a  quotient of the standard bicovariant calculus of $\OO_q(SU_{n})$.

\subsection{Quotients of the Standard Bicovariant Calculus on $\OO_q(SU_{n})$}

In this subsection we construct a left $\OO_q(SU_{n})$-covariant, right $\OO_q(L_m)$-covariant first-order calculus over $\OO_q(SU_{n})$. This calculus will be an essential tool used throughout the paper.

\begin{defn}
We denote by $\Omega^1_{q,\text{bic}}(SU_{n})$ the bicovariant first-order calculus induced on $\OO_q(SU_{n})$ by its coquasitriangular structure, and by $\Lambda^1_{q,\mathrm{bic}}$ the space of invariant forms. Moreover, we denote by 
$
\{b_{ij} \,|\, i,j = 1, \dots, n\}
$
the basis of $\Lambda^1_{q,\mathrm{bic}}$ dual to the functionals $Q_{ij}$.
\end{defn}

The restriction of the standard bicovariant calculus $\Om^1_{q,\text{bic}}(SU_{n})$ to the quantum Grassmannians can be shown to have  non-classical dimension even in the simplest cases of $\OO_q(\mathbb{CP}^{n})$ and $\OO_q(\textrm{Gr}_{4,2})$. Hence it cannot be isomorphic to the Heckenberger--Kolb calculus. We circumvent this problem by constructing a family of quotients of $\Om^1_{q,\text{bic}}(SU_{n})$. The proof of the following proposition is a simple and direct generalisation of \cite[Lemma 4.3]{MMF1} and so we omit it.

\begin{prop}
It holds that:
\begin{enumerate}
\item a  right submodule of $\Lambda^1_{q,bic}(SU_{n})$ is given by 
\begin{align*}
V_{m} := \spn_{\mathbb{C}}\{b_{ij} \,|\, i,j  = m+ 1, \ldots, n\},
\end{align*}
\item with respect to the right coaction
$$
\Delta_R := (\id \otimes \pi) \circ \mathrm{Ad_R}: \Lambda^1_{q,\mathrm{bic}}(SU_{n}) \to \Lambda^1_{q,\mathrm{bic}}(SU_{n}) \otimes \OO_q(L_m),
$$
the subspace $V_m$ is a right coideal of $\Lambda^1_{q,\mathrm{bic}}(SU_{n})$.
\end{enumerate}
\end{prop}

The following corollary now follows from a direct application of Takeuchi's equivalence, generalising \cite[Corollary 4.4]{MMF1}. 

\begin{cor} \label{cor:quotFODC}
A right $\OO_q(L_m)$-covariant, left $\OO_q(SU_{n})$-covariant, first-order differential calculus over $\OO_q(SU_{n})$ is given by 
\begin{align*}
\Om^1_q(SU_{n},m) := \Om^1_{q,bic}(SU_{n})/B_m,
\end{align*}
where we have denoted 
$$
B_m : = \unit \inv(A \otimes V_m) \sseq \Om^1_{q,bic}(SU_{n}).
$$
\end{cor}


\subsection{The Action of $Q$ on the Generators of $\OO_q(S^{n,m})$ and $\OO_q(\mathrm{Gr}_{n,m})$}

In this subsection we establish necessary conditions  for non-vanishing of the functions $Q_{ij}$ on the generators of $\OO_q(S^{n,m})$ and $\OO_q(\mathrm{Gr}_{n,m})$. These results are used in the next subsection to describe the restriction of the calculus $\Om^1_q(SU_{n},m)$ to the quantum Grassmannians.

\begin{lem}\label{lem:qcs}
~~~~  
\begin{enumerate}
\item For $I,J \sseq \{1, \dots, n\}$, we have $Q_{ij}([I|J]) \neq 0$  only if $I = J_{ij}$,  or equivalently only if  $J = I_{ji}$.
\item It holds that $Q_{ii}(z) = Q_{ii}(\overline{z})\inv  \neq 0$. Moreover $Q_{ii}=Q_{i'i'}$ for every $i,i'\in \{1,...,m\}.$
\end{enumerate}
\end{lem}
\begin{proof}
~~~
\begin{enumerate}
\item 
By Goodearl's formulae 
\begin{align*}
Q_{ij}([I|J]) =     \sum_{a=1}^{n} \sum_{A} r(u^i_a \oby [I|A])r([A|J] \oby u^a_j)
                       =   \sum_{a=1}^{n} r(u^i_a \oby [I|J_{aj}])r([J_{aj}|J] \oby u^a_j).
\end{align*}
Now $r(u^i_a \oby [I|J_{aj}])$ gives a non-zero answer only if $a \in I \cap J$ and $I_{ai} = J_{aj}$. This can happen only if  $I = J_{ij}$, or equivalently only if  $J = I_{ji}$. Hence $Q_{ij}([I|J])$ gives a non-zero answer only in the stated cases. 

\item It follows from Goodearl's formulae that, for $i=1, \dots, n$, 
\begin{align} \label{thuas}
Q_{ii}(z) = \sum_{a=1}^{n} r(u^i_a \oby [M|K])r(z_K \oby u^a_i) = r(u^i_i \oby z)r(z \oby u^i_i) \neq 0, 
\end{align}
as claimed. Moreover it follows from \cite[Lemma 2.1]{Goodearl2005} that $Q_{ii}(z)=Q_{i'i'}(z)$ for every $i,i'\in \{1,...,m\}.$
\end{enumerate}

\end{proof}

Building on this lemma, we next produce a set of necessary requirements for non-vanishing of the maps $Q_{ij}$ on the generators of $\OO_q(\mathrm{Gr}_{n,m})$. 

\begin{lem}\label{lem:QZIJ}
For $i,j = 1, \dots, n$, and $(i,j) \notin M^\perp \by M^\perp$, it holds that $Q_{ij}(z^{IJ}) \neq 0$ only in the three following cases:
\begin{enumerate}
 \item $I = M_{ij} \text{ and } J = M^\perp,$~~~~~~~~~~~~~~~~~~~~\,~ $ \text{  with } (i,j) \in M \by M^\perp,$
 \item $I = M \text{ ~  and } J = M^\perp_{ij},$~ ~~~~~~~~~~~\,\,~~~~~~~ $ \text{  with } (i,j) \in M^\perp \by M,$
\item $I=M$ ~ and $J = M^\perp$, ~ ~~~~~~~~~~~~~~~~~~  \text{  with } $i=j$.
\end{enumerate}
\end{lem}
\begin{proof}
Note first that 
\begin{align*}
Q_{ij}(z^{IJ}) = &\,  \sum_{|K| =m} \, \sum_{|L|=n-m}\, \sum_{a=1}^n r(u^i_a \oby [I|K][J|L])r(z^{KL} \oby u^a_j) \\
= &\,  \sum_{|K| = m} \, \sum_{|L|=n-m}  \, \sum_{a,b,c = 1}^n r(u^i_a \oby [J|L])r(u^a_b \oby [I|K])r(z_K \oby u^b_c)r(\overline{z}_L \oby u^c_j)\\
= & \, \sum_{|L|=n-m} \sum_{a, c=1}^n r(u^i_a \oby [J|L])Q_{ac}(z_I)r(\overline{z}_L \oby u^c_j.)
\end{align*}
Now $r(\overline{z}_L \oby u^c_j) \neq 0$ only if $c \leq j$. Thus $c \in M^\perp$ implies that $j \in M^\perp$. Now if  $c \neq j$, then $M^\perp_{cj}$ would contain $j$ as repeated element, and so, could not be equal to $L$ as required by Goodearl's formula. Hence
\begin{align*}
Q_{ij}(z^{IJ}) 
  = \sum_{a=1}^i  r(u^i_a \oby z_J)Q_{aj}(z_I)r(\overline{z}_L \oby u^j_j).
\end{align*}
If we assume that $I \neq M$, then by Lemma \eqref{lem:qcs},  $Q_{aj}(z_I) \neq 0$ only if \mbox{$(a,j) \in M \by M^\perp$} and $I = M_{aj}$. Since we have assumed that $(i,j) \notin M^\perp \by M^\perp$, we must have  that $i \in M$, and so,  $r(u^i_a \oby \overline{z}_J) \neq 0$ only  if $i=a$. Thus, if $I \neq M$, then  $Q_{ij}(z^{IJ}) \neq 0$ only if $I=M_{ij}$ and $J = M^\perp$, which gives us the first case. If we instead assume that $I =M$, then an analogous argument will shows  that $Q_{ij}(z^{IJ}) \neq 0$ only when $J = M^\perp_{ij}$, giving us the second and third cases.
\end{proof}
\begin{lem} It holds that
$Q_{ij}(z^{MM^\perp})  = Q_{ij}(1)$.
\end{lem}
\begin{proof}
We can derive this identity from Laplacian expansion as follows. In (\ref{Laplace}), setting $I =J = \{1, \dots, n\}$ and $J_1 = M$  gives 
\begin{align*}
1 = \mathrm{det}_n = \sum_{|I_1| = m} (-q)^{\ell(I_1,I_1^\perp)} z_M^{I_1}\overline{z}_{I_1^\perp} = \sum_{|I_1| = m} (-q)^{\ell(I_1,I_1^\perp)} \overline{z}_{I_1} \overline{z}_{I_1^\perp} = \sum_{|I_1| = r} (-q)^{\ell(I_1,I_1^\perp)}    z^{I_1I_1^\perp}.  
\end{align*}
Since $Q_{ij}(z^{IJ})$ is non-zero only in one of the three cases given above,  $Q_{ij}(z^{I_1I_1^\perp}) = 0$ unless $I_1 = M$. Hence, 
\begin{align*}
Q_{ij}(1) = \sum_{|I_1| =m} (-q)^{\ell(I_1,I_1^\perp)} Q_{ij}(z^{I_1I_1^\perp}) = Q_{ij}(z^{MM^\perp})
\end{align*}
as required.
\end{proof}

\subsection{A Restriction Presentation of the Heckenberger--Kolb Calculus} \label{section:HKRestriction}

Using the results of the previous subsection, we present  the Heckenberger--Kolb calculus as the restriction of the calculus $\Omega^1_q(SU_{n},m)$, reproduce its decomposition into $\Omega^{(1,0)}$ and $\Omega^{(0,1)}$, and introduce a basis of $\Phi\big(\Om^1_q(\mathrm{Gr})\big)$.

\begin{lem} \label{lem:vpmsubspaces} Denoting  
\begin{align*}
V\hol  := \{b_{ij} \,|\,   (i,j) \in M^\perp \by M \, \}, & & V\ahol  := \{b_{ij} \,|\,   (i,j) \in M \by M^\perp \, \},
\end{align*}
it holds that 
\begin{enumerate}
\item $V\hol$ and $V\ahol$ are non-isomorphic right $\OO_q(L_m)$-subcomodules of  $\Lambda^1$, 
\item $V\hol$ and $V\ahol$ are right $\OO_q(SU_{n})$-submodules of $\Lambda^1$.
\end{enumerate}
\end{lem}
\begin{proof}
The fact that $V\hol$ and $V\ahol$ are subcomodules of $\Lambda^1$ can be demonstrated using the explicit coset representatives of $b_{ij}$ presented in \cite[Lemma 4.2]{MMF1}: For example, for $(i,j) \in M \by M^\perp$, 
\begin{align*}
(\id \oby \pi) \circ \mathrm{Ad}_R[u^i_j] = & ~\, \sum_{a,b=1}^{n} [u^a_b] \oby \pi_{{n},m} \big(u^i_aS(u^b_j)\big)\\
= & \,\,\, \sum_{a,b=1}^{n} \sum_{x = 1}^m \sum_{y=m+1}^{n} [u^a_b] \oby \alpha_m\big(u^i_xS(u^y_j)\big) \oby \beta_m\big(u^x_aS(u^b_y)\big) \\
= & \,\,\, \sum_{a, x = 1}^m \sum_{b,y=m+1}^{n} [u^a_b] \oby \alpha_m\big(u^i_xS(u^y_j)\big) \oby \beta_m\big(u^x_aS(u^b_y)\big). 
 \end{align*}
Thus we see that $V^{(0,1)}$ is coinvariant under the right $\OO_q(L_m)$-coaction $(\id \oby \pi) \circ \mathrm{Ad}_R$. The fact that $V^{(1,0)}$ is also coinvariant is established similarly.

 Consider now the Hopf algebra map $\tau: \OO_q(L_m) \to \OO(U_1)$ uniquely determined by 
\begin{align*}
\tau(u^i_j\oby 1) :=  \delta_{ij} t,  & &   \tau(1\oby u^i_j) :=   \delta_{ij}.
\end{align*}

Composing the left $\OO_q(L_m)$-coactions on $V^{(1,0)}$ and $V^{(0,1)}$ with $\id \otimes \tau$ we get a right $\OO(U_1)$-coaction on both spaces which we denote by $\Delta_{R,\OO(U_1)}$. From the calculation of the $\OO_q(L_m)$-coactions above, it follows that $\Delta_{R,\OO(U_1)}$ acts according to 
\begin{align*}
\Delta_{R,\OO(U_1)}[u^i_j] =[u^i_j] \otimes t \quad \textrm{for} \; [u^i_j]  \in V^{(0,1)}
\end{align*}
and
\begin{align*}
\Delta_{R,\OO(U_1)}[u^i_j] = [u^i_j] \otimes t^{-1} \quad \textrm{for} \; [u^i_j]  \in V^{(1,0)}.
\end{align*}
Thus we see that, since the two comodules are not isomorphic as $\OO(U_1)$-modules, they cannot be isomorphic as $\OO_q(L_m)$-modules. 

Finally, the fact that $V\hol$ is an $\OO_q(SU_{n})$-submodule follows from the calculation
\begin{align*}
Q_{ab}(u^i_ju^k_l) = 0, & & \text{ for all } (i,j) \in M^\perp \by M; \, (a,b) \notin M \by M^\perp, \; k,l = 1, \dots, n.
\end{align*}
A similar argument shows that $V\ahol$ is a submodule.
\end{proof}

\begin{prop}\label{prop:resthk} The restriction of  $\Om^1_q(SU_{n},m)$ to $\OO_q(\mathrm{Gr}_{n,m})$ is the Heckenberger--Kolb calculus, with the decomposition into subcalculi given by 
\begin{align*}
\Om\hol := \OO_q(SU_{n}) \square_{\mathcal{O}_q(L_m)} V\hol, & &  \Om\ahol := \OO_q(SU_{n}) \square_{\mathcal{O}_q(L_m)} V\ahol,
\end{align*} 
\end{prop}
\begin{proof}
Denote by $\Om^1_q(\mathrm{Gr}_{n,m})$ the restriction of $\Om^1_q(SU_{n},m)$ to $\OO_q(\mathrm{Gr}_{n,m})$.  It follows from Proposition \ref{prop:PPQPB} and Lemma \ref{lem:vpmsubspaces} that the image of $\Phi(\Om^1(\mathrm{Gr}_{n,m}))$ in $\Lambda^1$ is given by the linear span of the elements of the form $[b^+]$, for $b \in \OO_q(\mathrm{Gr}_{n,m})$. Consider next the right $\OO_q(\mathrm{Gr}_{n,m})$-ideal generated by the elements 
$(z^{IJ})^+$, for $|I|= m, |J| = n-m$. Noting that the quotient of $\OO_q(\mathrm{Gr}_{n,m})$ by this ideal is one-dimensional, we see that $\OO_q(\mathrm{Gr}_{n,m})^+$ is generated as a right $\OO_q(\mathrm{Gr}_{n,m})$-module by the the elements $(z^{IJ})^+$. Lemma \ref{lem:QZIJ} implies that
\begin{align*}
\spn_{\mathbb{C}}\big\{[(z^{IJ})^+] \, |\,  |I| = m, |J|=n-m \big\} \sseq V\hol \oplus V\ahol.
\end{align*}
Thus, since Lemma \ref{lem:vpmsubspaces} above tells us that $V\hol \oplus V\ahol$ is a right submodule of $\Lambda^1$,  we must have that 
\begin{align*}
\big\{[b^+] \, |\, b \in \OO_q(\mathrm{Gr}_{n,m})\big\} \sseq  V\hol \oplus V\ahol.
\end{align*}
Moreover, Lemma \ref{lem:QZIJ} tells us that each non-zero $[\big(z^{IJ}\big)^+]$ is contained in either $V\hol$ or $V\ahol$. Thus, since $V\hol$ and $V\ahol$  are each submodules of $V\hol \oplus V\ahol$, we must have a decomposition
\begin{align*}
\Phi\left(\Om^1(\mathrm{Gr}_{n,m})\right) = W\hol \oplus W\ahol :=  \big(\Phi(\Om^1(\mathrm{Gr}_{n,m}) \cap V\hol\big) \oplus  \big(\Phi(\Om^1(\mathrm{Gr}_{n,m}) \cap V\ahol\big)
\end{align*}
as well as a corresponding decomposition of calculi. Since Lemma \ref{lem:vpmsubspaces} tells us that $V\hol$ and $V\ahol$ are non-isomorphic  as left $\OO_q(L_m)$-comodules, $\Psi(V\hol)$ and $\Psi(V\ahol)$ are non-isomorphic as objects in ${}^{\OO_q(SU_n)}_{\OO_q(\mathrm{Gr}_{n,m})}\mathrm{Mod}$, and so, the two calculi cannot be isomorphic.

The classification of calculi over the quantum Grassmannians given in Theorem \ref{thm:HKClass} tells us that there can exist no non-trivial calculus of dimension strictly less than $\dim\big(\Om^1_q(\mathrm{Gr}_{n,m})\big)$. Thus it follows from the inequality
\begin{align*}
\dim\big(W\hol\big)  \leq \dim\!\big(V\hol\big) = m(n-m) = \frac{1}{2} \dim\big(\Om^1_q(\mathrm{Gr}_{n,m})\big)
\end{align*}
that $W\hol = V\hol$,  or $W\hol = 0$. That the former is true follows from the fact that the coset  $[z^{MM^{\perp}_{m+1,1}}]$ is non-zero. This can be seen explicitly from the fact that  $Q_{m+1,1}\big(z^{MM^\perp_{m+1,1}}\big)$ is equal to 
\begin{align*}
& \sum_{a=1}^{n}\sum_{b=1}^{n}\sum_{|A|=m}\sum_{|B|=n-m}r(u^{m+1}_a \oby [M^\perp_{m+1,1}|B])r(u^{a}_b \oby [M|A])r  (z_A\oby u^{b}_c ) r(\overline{z}_B\oby u^{c}_1),
\end{align*}
which, since the last factor is non-zero only for $c=1$, is equal to 
\begin{align*}
 r(u^{m+1}_1 \oby \overline{z}_{M^\perp_{m+1,1}}) r\left(u^{1}_1 \oby \overline{z}\right)r  (z\oby u^{1}_1 ) r(\overline{z}\oby u^{1}_1 ) 
\neq   0.
\end{align*}
A similar argument establishes that $W\ahol = V\ahol$.
Since the associated calculi are non-isomorphic,  their direct sum must be isomorphic to the Heckenberger--Kolb calculus. \end{proof}


\section{A Borel--Weil Theorem for the Quantum Grassmannians} \label{section:BorelWeil}

In this section we prove the main result of the paper, namely the Borel-Weil theorem for the quantum Grassmannians. To do so we realize the unique $\OO_q(SU_n)$-covariant connection on the Heckenberger-Kolb calculus as a principal connection acting on a quantum principal bundle over the Grassmannian.
\subsection{A Quantum Principal Bundle}
In this subsection we prove that the sphere bundle $\mathcal{O}_q(S^{n,m})$ together with the restriction of the first order differential calculus $\Omega^1_q(SU_n),$ denoted as $\Omega^1_q(S^{n,m}),$ gives a quantum principal bundle. Then we prove that the universal principal connection $\Pi_i$ descends to a principal connection on $\Omega^1_q(SU_n).$

\begin{defn}
We denote by $\Omega^1_q(S^{n,m})$ the restriction of $\Omega^1_q(SU_n)$ to a first-order differential calculus on $\OO_q(S^{n,m})$. 
\end{defn}

We observe that since $\Omega^1_q(SU_n,m)$ restricts to the first-order Heckenberger--Kolb calculus $\Omega^1_q(\mathrm{Gr}_{n,m})$ on $\OO_q(\mathrm{Gr}_{n,m})$, then $\Omega^1_q(S^{n,m})$ must also restrict to $\Omega^1_q(\mathrm{Gr}_{n,m})$ on $\OO_q(\mathrm{Gr}_{n,m})$.

\begin{prop}
The first-order differential calculus $\Omega^1_q(S^{n,m})$ is left $\OO_q(SU_n)$-covariant and right $\OO_q(L_m)$-covariant. Hence a quantum principal bundle is given by the pair  
\begin{align} \label{eqn:thebundle}
\left(\OO_q(S^{n,m}),\Omega^1_q(S^{n,m})\right).
\end{align}
\end{prop}
\begin{proof}
The left and right covariance of $\Omega^1_q(S^{n,m})$ is inherited directed from the left and right covariance of $\Omega^1_q(SU_n,m)$. Hence by Proposition \ref{prop:PPQPB} the given pair is a quantum principal bundle.
\end{proof}
    
We will now prove that the universal principal connection $\Pi_i$ descends to a principal connection on $\Omega_q^1(S^{n,m})$. We begin with a general technical lemma for principal pairs.
\begin{lem}\label{lem:jnothom}
Denote $J_{\mathrm{ver}} := J\cap i(T^+)$. If $J$ is not homogeneous with respect to the decomposition given in \eqref{eqn:Pdecomposition} then 
\begin{align}\label{eqn:dimbound}
\dim(\Lambda^1_H) <  \dim\!\left(T^+\!/\pi(J_{\mathrm{ver}})\right)\!.
\end{align}
\end{lem}
\begin{proof}
We can decompose any element $j\in J$ as $j=j_0+j_1$, with $j_0 \in \mathrm{ker}(\pi)$ and $j_1 \in i(T^+).$ Now we have
\begin{align}\label{eqn: rhoj}
    \pi(j)=\pi(j_1),
\end{align}
for every $j \in J.$ If $J$ is not homogeneous, there exist $j\in J$ such that $j_1 \notin J_{\mathrm{ver}}$ and, since $i:T^+\rightarrow P^+$ is injective, equation \eqref{eqn: rhoj} implies
\begin{align*}
     \dim \left(T^+/\pi(J)\right)     <  \dim \left(T^+\!/\pi(J_{\mathrm{ver}})\right),
    \end{align*}
as claimed.
    \end{proof}

\begin{prop} ~~~~
\begin{enumerate}
    \item  The universal principal connection $\Pi_i$ descends to a left $\OO_q(SU_n)$-principal connection for the bundle given in  \eqref{eqn:thebundle}.
    \item The associated space of vertical forms $V_q(S^{n,m})_{\mathrm{ver}}$ is $1$-dimensional, spanned by the coset $[z-1]$, or equivalently by the coset $[\overline{z}-1]$, and it is trivial as a $\OO_q(L_m)$-comodule. 
    \item The connection $\Pi_i$ is the unique left $\OO_q(SU_n)$-covariant principal connection for the bundle.
\end{enumerate}
\end{prop}
\begin{proof}
Denoting $\lambda:=Q_{ii}(\overline{z}),$
it follows from \cite[Lemma 2.1]{Goodearl2005} that $\lambda \neq 1$. Thus we see  that $[\overline{z} - 1] \neq 0$. Moreover, since we have that 
$$
[\pi(\overline{z} - 1)] = [t-1] \neq 0,
$$
implying that $\Lambda^1_H$ is a non-trivial space spanned by the coset $[t-1]$.
\\To prove that $J$ is homogeneous with respect to the decomposition given by \eqref{eqn:Pdecomposition}, let  $J_{\mathrm{ver}}=J\cap i(T^+).$ The elements
$$ 
\left(\overline{z}-\lambda\right)\!\left(\overline{z} - 1\right)\!\overline{z}^k\!,~ \left(z -\lambda\right)\!\left(z-1\right)\!z^k
$$
and 
$$ 
\left(\overline{z}-1\right)-\frac{\lambda-1}{\lambda^{-1}-1}\left(z-1\right)
$$
are contained in $J_{\mathrm{ver}}$. Thus we can conclude that
\begin{align*}
    \dim\!\left(T^+\!/\pi (J_{\mathrm{ver}})\right)\leq 1.
\end{align*}
If $J$ were not homogeneous with respect to the decomposition Lemma \eqref{lem:jnothom} would then imply that $T^+\!/\pi (J)$ was trivial, but since we know that $\pi(\overline{z}-1)\neq 0$, we must conclude that $J$ is homogeneous.

Recall that $V^{(0,1)}$ and $V^{(1,0)}$ are irreducible and non-isomorphic, and for all but the very special case of the Podle\'s sphere, have dimension strictly greater than $1$. Hence neither $V^{(1,0)}$ nor $V^{(0,1)}$ can be be trivial as a $\OO_q(L_m)$-comodule. For the special case of the Podle\'s sphere, both $V^{(1,0)}$ and $V^{(0,1)}$ are $1$-dimensional, but are well-known to be non-trivial (see for example \cite[\textsection 3]{Maj}.) All this means that that $V_q(S^{n,m})$ decomposes into three distinct irreducible sub-comodules, and hence that there can exists precisely one left $\OO_q(L_m)$-comodule map with kernel $V^{1}_q(\mathrm{Gr}_{n,m})$, and hence precisely one left $\OO_q(SU_n)$-covariant principal connection for the bundle.
\end{proof}

Following the general discussion in \textsection  \ref{subsection:ConnFromPrinCs}, the principal connection $\Pi:\Omega^1_q(S^{n,m}) \mapsto \Omega^1_q(S^{n,m})$ gives a left $\OO_q(SU_n)$-covariant connection 
\begin{align} \label{eq:GrassPrinConn}
\nabla: \EE_k \mapsto \Omega^1_q(\mathrm{Gr}_{n,m}) \otimes \EE_k,
\end{align}
for each $k \in \mathbb{Z}$. The following corollary is now implied by the uniqueness results of connections identified in Theorem \ref{thm:covconnection}.

\begin{cor}
For every line bundle $\EE_k$, the connection $\nabla$ coincides with the connection identified in Theorem \ref{thm:covconnection}. Moreover, the associated $(0,1)$-connection is the holomorphic structure identified in Theorem \ref{thm:covconnection}.
\end{cor}

\subsection{The Borel--Weil Theorem} \label{subsection:BW}

In this subsection we prove the Borel-Weil theorem for quantum Grassmannians. We begin with a short technical lemma used at a number of points in the proof of the main theorem. In the lemma the projection 
$$
\proj^{(0,1)}: \Omega^1_q(S^{n,m}) \mapsto \OO_q(S^{n,m})\Omega^{(0,1)}\OO_q(S^{n,m})
$$
associated to the decomposition 
$$
\Omega^1_q(S^{n,m}) = \OO_q(S^{n,m})\Omega^{(1,0)}\OO_q(S^{n,m}) \oplus \OO_q(S^{n,m})\Omega^{(0,1)}\OO_q(S^{n,m}) \oplus \Omega^{1}_q(S^{n,m})_{\mathrm{ver}}
$$
is considered.

\begin{lem} \label{lem:projrightmodpve}
It holds that 
\begin{enumerate}
    \item $
\sigma\Big(\Phi_P(\OO_q(S^{n,m})\Omega^{(0,1)}\OO_q(S^{n,m}))\Big) = \sigma\Big(\Phi_P(\OO_q(S^{n,m})\Omega^{(0,1)}\Big) = V^{(1,0)},
$
\item $
\sigma\Big(\Phi_P(\OO_q(S^{n,m})\Omega^{(1,0)}\OO_q(S^{n,m}))\Big) = \sigma\Big(\Phi_P(\OO_q(S^{n,m})\Omega^{(1,0)}\Big) = V^{(0,1)}.
$
\item The projection $\mathrm{proj}^{(0,1)}$ is a right $S_q(\mathrm{Gr}_{n,m})$-module map.
\end{enumerate}
\end{lem}
\begin{proof}
It follows from the presentation of $\Omega^{(1,0)}$ and $\Omega^{(0,1)}$ given in Proposition \ref{prop:resthk} that 
\begin{align*}
\sigma([\omega]) \in V^{(1,0)}, & & \mu([\omega]) \in V^{(0,1)}, & & ~~ \omega \in \Omega^{(1,0)}, \, \nu \in \Omega^{(0,1)}.      
\end{align*}
Thus we see that the equalities in 1 and 2  follow from the fact that $V^{(1,0)}$ and $V^{(0,1)}$ are right $\OO_q(SU_n)$-submodules of $V_q(S^{n,m})$, as established in Lemma \ref{lem:vpmsubspaces}. Lemma \ref{lem:vpmsubspaces} also tells us that $\mathrm{proj}^{(0,1)}$ is a right $S_q(\mathrm{Gr}_{n,m})$-module map if 
\begin{align*}
V^0S_q(\mathrm{Gr}_{n,m}) \sseq V^{(1,0)} \oplus V^0.
\end{align*}
That this is true follows from the fact that 
$$
Q_{ij}\left(\big(\overline{z} - 1\big)\overline{z}_I\right) = 0,
$$
for all $(i,j) \in M \times M^{\perp}$.
\end{proof}

\begin{thm}[Quantum Grassmannian Borel--Weil] \label{thm:BW}
For all $k \in \mathbb{Z}_{>0}$, it holds that 
\begin{align*}
1. ~ H^0_{\adel}(\mathcal{E}_k) = V_{k\varpi_{n-m}}, & & 2. ~ H^0_{\adel}(\mathcal{E}_{-k}) = 0,
\end{align*}
\end{thm}
\begin{proof}
~~~

1. ~~  Recalling that $z = [M^{\perp}|M^{\perp}]$ we see that
\begin{align*}
(\id \otimes \sigma)\unit\!\left(\exd z\right) 
=  \sum_{|A|= n -m} ~ \sum_{(i,j) \in M\times M^\perp}  [M^{\perp}|A] \otimes  Q_{ij}([A | M^{\perp}])b_{ij}.
\end{align*}
It follows from Lemma \ref{lem:qcs} that $Q_{ij}([A | M^{\perp}]) = 0$, for all $(i,j) \notin M^{\perp} \times M$. The results of \textsection \ref{section:HKRestriction} now imply that 
\begin{align*}
\exd \overline{z} \in \Omega^{(1,0)}_q(S^{n,m})_{\mathrm{hor}} \oplus \Omega^1_q(S^{n,m})_{\mathrm{ver}} = \ker(\mathrm{proj}^{(0,1)}).
\end{align*}
Note next that, for any $k \in \mathbb{Z}_{>0}$, Proposition 1.34 implies  
\begin{align*}
\adel_{\EE_{k}}(\overline{z}^k) =  \mathrm{proj}^{(0,1)} \circ \exd(\overline{z}^k) = \sum_{a+b = k-1~} \overline{z}^a \, \mathrm{proj}^{(0,1)}(\exd \overline{z}) \overline{z}^b = 0.
\end{align*}
Thus all powers of $\overline{z}$ are holomorphic. Recalling that  $\adel_{\EE_k}$ is a left $\OO_q(SU_n)$-comodule map, we see that the inclusion of $V_{k\varpi_{n-m}}$ in  $H^0_{\adel}(\mathcal{E}_k)$ now follows from Schur's lemma.

To show the opposite inclusion, we need to consider the decomposition of $\EE_k$ into irreducible left $\OO_q(SU_n)$-subcomodules. As discussed in Appendix \ref{Appendix:GensofGrassSphere}, each irreducible subcomodule contains an element of the form $b\overline{z}^k$, for some $b \in \OO_q(\mathrm{Gr}_{n,m})$. Let us assume that one of these elements $b\overline{z}^k$, for $b \notin \mathbb{C}1$, were holomorphic. 
This would imply that
     \begin{align*}
         \adel_{\EE_k}\!(b \overline{z}^{\,k}) =  \adel b \otimes \overline{z}^{\,k} + b \otimes \adel_{\EE_k}(\overline{z}^{\,k})
         =  \adel b \otimes \overline{z}^{\,k}.
     \end{align*}
Noting  that $\Omega^1_q\!\left(S^{n,m}\right)$ is a torsion-free right $\OO_q(\mathrm{Gr}_{n,m})$-module,  that $\mathcal{E}_k$ is projective as left $\OO_q(\mathrm{Gr}_{n,m})$-module and that $\OO_q(\mathrm{Gr}_{n,m})$ has no zero-divisors, we see that Liouville's theorem implies that 
$$
\adel b \otimes \overline{z}^{\,k} \neq 0.
$$
Thus we see that no such holomorphic element exists, and that by Schur's lemma the first claimed identity holds.

\bigskip

2. ~~ We now come to the line bundles $\mathcal{E}_{-k}$, for $k \in \mathbb{Z}_{>0}$. Assume there exists a non-zero holomorphic element $e \in \mathcal{E}_{-k}$. 
Then for any $|I| =n-m$, we would have
\begin{align*}
~~~~ \adel\big(e\overline{z}_I^{\,k}\big) = & \,  \mathrm{proj}^{(0,1)} \Big(\exd\big(e\overline{z}_I^{\,k}\big) \Big)
=  \, \mathrm{proj}^{(0,1)}\!\left((\exd e)\overline{z}_I^{\,k} + e\exd(\overline{z}_I^{\,k})\right)
=  0,
\end{align*}
which is to say 
$$
e\overline{z}_I^k \in H^0_{\adel,q}(\mathrm{Gr}_{n,m}).
$$ 
Lemma 1.23 implies that $e \overline{z}_I^{\,k}$ must be a non-zero scalar multiple of $1$. Denoting
$$
e \overline{z}_I^{\,k}  =: \mu_I 1,
$$

we see that $\mu_I^{-1} \overline{z}_I^{\,k}$ is a right inverse of $e$, for any $I$. However, since $\OO_q(SU_n)$ has no zero divisors, right inverses are unique. Thus to avoid contradiction we are forced to conclude that $\EE_{-k}$ contains no holomorphic elements. 
\end{proof}

\subsection{Opposite Complex Structure}

In this subsection we consider the opposite complex structure $\ol{\Om}^{(\bullet,\bullet)}$ of the complex structure $\Omega^{(\bullet,\bullet)}$. Composing the connection given in \eqref{eq:GrassPrinConn} with projection onto $\Omega^{(1,0)} \otimes_{\OO_q(\mathrm{Gr}_{n,m})} \EE_k$, gives a left $\OO_q(SU_n)$-covariant flat $(1,0)$-connection 
$$
\del_{\EE_k}: \EE_k \to \Omega^{(1,0)} \otimes_{\OO_q(\mathrm{Gr}_{n,m})} \EE_k,
$$
which is to say a covariant holomorphic structure with respect to $\Omega^{(\bullet,\bullet)}$ (see \cite[Remark 4.8]{HVBQFM} for more details). As a carefully reading of \textsection \ref{subsection:BW} will confirm, the proof of the Borel--Weil theorem for the complex structure $\Om^{(\bullet,\bullet)}$ carries over directly to the opposite setting. 

\begin{thm} \label{thm:BW}
For the holomorphic line bundle $(\EE_k, \adel_{\EE_k})$, it holds that 
\begin{align} 1. ~ H^0_{\del}(\mathcal{E}_k) = 0, & & 
2. ~ H^0_{\del}(\mathcal{E}_{-k}) \simeq V_{k\varpi_m}, 
\end{align}
for all $k \in \mathbb{Z}_{>0}$.
\end{thm}

This result will be used in \textsection \ref{section:Kahler} when we discuss the implications of the Borel--Weil theorem for the noncommutative K\"ahler geometry of the quantum Grassmannians.

\subsection{The Twisted Homogeneous Coordinate Ring}

In classical complex geometry any positive line bundle $\EE$ over a compact K\"ahler manifold $M$ gives an embedding of $M$ into complex projective space. Moreover, the homogeneous coordinate ring  $S(M)$ associated to this embedding is isomorphic as a graded algebra to
\begin{align} \label{eqn:HoloHCR}
\bigoplus_{k \in \mathbb{Z}_{\geq 0}} H^0_{\adel}(\EE^{\otimes k}).
\end{align}

For the Grassmannians the line bundle $\EE_1$ is well-known to be positive and the associated projective embedding is the Pl\"ucker embedding. Thus for the quantum Grassmannians it is natural to ask if the ring 
\begin{align} 
\bigoplus_{k \in \mathbb{Z}_{\geq 0}} H^0_{\adel}(\EE_{k})
\end{align}
can be identified with the twisted homogeneous coordinate ring $S_q(\mathrm{Gr}_{n,m})$. The following proposition shows that this is indeed the case and that the result is a direct implication of the quantum Borel--Weil theorem. This generalises earlier work in \cite{KLVSCP1,KKCP2,KKCPN} realising the quantum projective plane in terms of the noncommutative complex geometry of quantum projective space.

\begin{prop} \label{prop:HHCR}
Considering each line bundle $\EE_k$ as a subspace of $\OO_q(SU_n)$, it holds that 
$$
\bigoplus_{k \in \mathbb{Z}_{>0}} H^0_{\adel}(\mathcal{E}_k) = S_q(\mathrm{Gr}_{n,m}).
$$
\end{prop}
\begin{proof}
It follows from the Borel--Weil theorem that
\begin{align*}
\bigoplus_{k \in \mathbb{Z}_{\geq 0}} H^0_{\adel}(\EE_k)  \sseq S_q(\mathrm{Gr}_{n,m}).
\end{align*}
Let us now show the opposite inclusion. Take $e \in H^0_{\adel}(\EE_k), \, e' \in H^0_{\adel}(\EE_l)$, and observe that 
\begin{align*}
 \adel_{\EE_{k+l}}(ee') = & \, (\id -\Pi)^{(0,1)}\!\left(\exd(ee')\right) \\
 = & \, (\id -\Pi)^{(0,1)}\!\left((\exd e)e' + e\exd(e')\right) \\
 = & \, (\id -\Pi)^{(0,1)}(\exd e)e' + e\Pi(\exd e') \\
 = & \, 0,
\end{align*}
where we have used the fact that $\Pi$ is a right $S_q(\mathrm{Gr}_{n,m})$-map, as established in Lemma \ref{lem:projrightmodpve}. Thus we see that the direct sum $\bigoplus_{k \in \bZ_{\geq 0}} H^0_{\adel}(\mathcal{E}_k)$ is a subalgebra of $\OO_q(SU_n)$. Now by the Borel--Weil theorem the generators $\overline{z}_I$ of $S_q(\mathrm{Gr}_{n,m})$ are contained in $H^0_{\adel}(\mathcal{E}_1)$. Thus we get the opposite inclusion
$$
S_q(S^{n,m}) \sseq \bigoplus_{k \in \mathbb{Z}_0} H^0_{\adel}(\mathcal{E}_k),
$$  
and hence equality of the two algebras.
\end{proof} 

\begin{remark} \label{rem:positiveEE1}
Proposition \ref{prop:HHCR} hints at the possible existence of a noncommutative generalisation of the classical GAGA correspondence. (See \cite[\textsection 7]{BS} for a detailed  speculation about how  such a general picture might look.) The task of finding a unified understanding of the noncommmutative complex and projective geometry of the quantum Grassmannians in terms of such a formal framework is a major undertaking. A tractable first step is the formulation of a noncommutative version of positivity for noncommutative holomorphic line bundles. This is the subject of ongoing research and is discussed in some detail in the next section. 
\end{remark}

\section{Applications to the K\"ahler Structure of the Quantum Grassmannians} \label{section:Kahler}

The Heckenberger--Kolb calculi of the quantum Grassmannians possess a rich noncommutative  generalisation of the K\"ahler geometry of the classical Grassmannians. Motivated by a desire to understand this noncommutative K\"ahler geometry, a general theory of noncommutative K\"ahler structures was introduced in \cite{MMF3} and further developed in \cite{DOKSS,HVBQFM,OSV}. By feeding the numerical invariants given by Theorem \ref{thm:BW} into this general theory one can establish a number of new results for the noncommutative K\"ahler geometry of the quantum Grassmannians. This section presents a brief summary of these applications, referring the interested reader to relevant papers as we go. In particular we recall the generalisation of positivity to the noncommutative setting and verify that $\EE_1$ is a positive (see Remark \ref{rem:positiveEE1} above).

\subsection{K\"ahler Structures}

We begin with the definition of a noncommutative K\"ahler structure. It abstracts the properties of the fundamental $(1,1)$-form of a K\"ahler metric \cite[\textsection 3.1]{HUY}. 

\begin{defn} An {\em Hermitian structure} $(\Om^{(\bullet,\bullet)}, \sigma)$ for a differential $*$-calculus $\Om^{\bullet}$ of even total degree $2n$  is a pair  consisting of  a complex structure  $\Om^{(\bullet,\bullet)}$ and  a central real $(1,1)$-form $\sigma$, 
such that, with respect to the {\em Lefschetz operator}
\begin{align*}
L:\Om^\bullet \to \Om^\bullet,  & &   \omega \mapsto \sigma \wedge \omega,
\end{align*}
isomorphisms are given by
\begin{align*}
L^{n-k}: \Om^{k} \to  \Om^{2n-k}, & & \text{ for all } k = 0, \dots, n-1.
\end{align*}

A {\em K\"ahler structure} for a differential $*$-calculus is an Hermitian structure $(\Om^{(\bullet,\bullet)},\kappa)$ such that $\kappa$ is closed, which is to say, $\exd \kappa = 0$. 
\end{defn} 

Associated to every K\"ahler structure we have a rich noncommutative generalisation of the classical structures of K\"ahler geometry. For example, one gets a direct noncommutative generalisation of Lefschetz decomposition allowing for the definition of a noncommutative  Hodge map $\ast_{\kappa}$ through a generalisation of the classical Weil formula. One can then use the Hodge map to define a noncommutative K\"ahler metric
\begin{align*}
    g: \Omega^{\bullet} \times \Omega^{\bullet} \to B:= \Omega^0, & & (\omega,\nu) \mapsto \ast_{\kappa}(\omega \wedge \ast_{\kappa}(\nu^*)).
\end{align*}
An Hermitian module over $B$ is a pair $(\F,h)$, where $\F$ finitely generated projective left $B$-module together with an appropriately defined sesquilinear map
$$
h: \F \times \F \to B.
$$
When such an $\F$ is additionally endowed with a holomorphic structure $\adel_{\F}$, Beggs and Majid showed that one can mimic the classical Chern connection construction and produce a connection 
$$
\nabla:\F \to \Omega^1(B) \otimes_B \F,
$$
extending the $(0,1)$-connection $\adel_{\F}$ (see \cite[\textsection 8.6]{BeggsMajid:Leabh} for details.) This gives a natural extension of the classical notions of positive and negative line bundles. Explicitly, we say that an Hermitian holomorphic vector bundle $\F$ over $B$ is \emph{positive}, or respectively \emph{negative}, if the Chern curvature satisfies
$$
\nabla^2(f) = - \mathbf{i}\theta \kappa \otimes f,
$$
where $\theta$ is an element of $\mathbb{R}_{>0}$, or respectively an element of $\mathbb{R}_{<0}$. Under additional positivity assumptions on $g$, and the existence of a state $\phi:B \to \mathbb{C}$, one can then  build upon this definition to prove a direct noncommutative generalisation of the classical Kodaira vanishing theorem \cite[\textsection 8]{OSV}.

\subsection{Positive Line Bundles and Higher Cohomologies}

As established in \cite{MarcoConj}, for all but a finite (possibly empty)  set of values of $q$, the Heckenberger--Kolb calculi over $\OO_q(\mathrm{Gr}_{n,m})$ posses a covariant K\"ahler structure that is unique up to real scalar multiple. Moreover, it was shown in \cite[Proposition 7.13]{DOSFred} that there exists an open interval $I$  around $1$ such that, for all $q \in I$, the associated K\"ahler metric $g$ is positive definite. Combining this with the fact that the two first-order Heckenberger--Kolb calculi are irreducible, one can conclude from \cite[Theorem 3.3]{HVBQFM} that each holomorphic line bundle over $\OO_q(\mathrm{Gr}_{n,m})$ is either positive, negative or flat. Indeed it follows from \cite[Corollary 3.4]{DOSFred} that, for all $l \in \mathbb{Z}$,  we have 
\begin{enumerate}
\item  if $H^0_{\adel}(\EE_k) \neq 0$ and $H^0_{\del}(\EE_k) = 0$, then $\EE_k$ is positive,
\item  if $H^0_{\adel}(\EE_k) = 0$ and $H^0_{\del}(\EE_k) \neq 0$, then $\EE_k$ is negative.
\end{enumerate}
Thus we see that the Borel--Weil theorem for $\OO_q(\mathrm{Gr}_{n,m})$  implies that, for all $k \in \mathbb{Z}_{>0}$, the bundle $\EE_k$ is
positive, while the bundle $\EE_{-k}$ is negative.

Finally, in \cite[Theorem 4.18]{HVBQFM} the noncommutative Kodaira vanishing theorem for positive line bundles was applied to $\EE_k$. In particular, vanishing of higher cohomologies was concluded:
\begin{align*}
    H^{(0,i)}_{\adel}(\EE_k) = 0, & & \textrm{ for all } i \in \mathbb{Z}_{>0}.
\end{align*}
Thus we see that for positive line bundles, the classical Bott--Borel--Weil theorem \cite{BottBW}  extends to the quantum setting.

\subsection{Square Integrable Forms, Fredholm Operators, and $C^1$-Line Modules}

Composing the Haar state of $\OO_q(SU_n)$ with the K\"ahler metric $g$, gives an inner product 
\begin{align*}
\langle \cdot,\cdot\rangle: \Omega^{\bullet} \times \Omega^{\bullet} \to \OO_q(\mathrm{Gr}_{n,m}), & & (\omega,\nu) \mapsto \haar\left(g(\omega,\nu^*)\right)\!.
\end{align*}
The anti-holomorphic derivative $\adel$ is adjointable with respect to the inner product. Denoting the adjoint by $\adel^{\dagger}$, a $q$-deformation of the Dolbeault--Dirac operator of $\mathrm{Gr}_{n,m}$ \mbox{is given by}
$$
D_{\adel} := \adel + \adel^{\dagger}.
$$
We denote by $L^2(\Omega^{\bullet})$ the Hilbert space completion of $\Omega^{\bullet}$ with respect to the inner product. The Dolbeault--Dirac operator $D_{\adel}$ is an essentially self-adjoint operator on $L^2(\Omega^{\bullet})$, and the commutators $[D_{\adel},b]$ are bounded, for all $b \in \OO_q(\mathrm{Gr}_{n,m})$. Moreover, $L^2(\Omega^{\bullet})$ carries an evident representation $\rho$ of $\OO_q(\mathrm{Gr}_{n,m})$ by bounded operators. This means that the triple $(D_{\adel},L^2(\Omega^{\bullet}),\rho)$ is a natural candidate for a spectral triple, Connes' $C^*$-algebraic notion of a noncommutative Riemannian manifold \cite{RennieSpecTrip,Connes}. The difficulty here lies in verifying the compact resolvent condition for $D_{\adel}$. This has been verified by direct calculation of the spectrum for the special case of quantum projective space, but the extension of this result to the quantum Grassmannians is technically prohibitive.

Approaching the problem from a different point of view, we can consider Dolbeault--Dirac operators twisted by holomorphic line bundles:
\begin{align} \label{eqn:twistedDD}
D_{\adel_{\EE_{k}}}: \Omega^{(0,\bullet)} \otimes_{\OO_q(\mathrm{Gr}_{n,m})} \EE_k \to \Omega^{(0,\bullet)} \otimes_{\OO_q(\mathrm{Gr}_{n,m})} \EE_k, & & \textrm{ for } k \in \mathbb{Z}
\end{align} 
Since we now know that $\EE_{-k}$ is a negative line bundle, it follows from the results of \cite[\textsection 6]{DOSFred} that the eigenvalues of the $\EE_{-k}$-twisted Dolbeault--Dirac operator are bounded below by a non-zero positive number. Hence it follows that the twisted Dolbeault--Dirac operator is Fredholm. While the Fredholm property is significantly weaker than having compact resolvent, the conceptual geometric proof of this analytic property is a significant step in the right direction.

We finish by considering the domain of the closure of the essentially self-adjoint operator defined in \eqref{eqn:twistedDD}. This is a direct $q$-deformation of the vector space of $C^1$-sections, and it is natural to ask if the kernel of $\adel_{\EE_{-k}}$ contains additional elements, which is to say, whether there exist additional holomorphic elements in the $C^1$-completion of $\EE_k$. However, this cannot be the case since the Hilbert space completion of $\EE_k$ admits an orthonormal basis consisting of eigenvectors of the twisted Dolbeault--Dirac operator, which is to say, the operator is diagonalisable \cite[Proposition 4.19]{DOSFred}. This generalises the result for the special case of quantum projective space established in  \cite{KLVSCP1,KKCP2, KKCPN}.


\appendix

\section{A Generating Set for the Quantum Grassmannians} \label{Appendix:GensofGrassSphere}

In this appendix we recall the alternative quantised enveloping algebra presentation of the quantum Grassmannians introduced in \cite{DijkStok}, and equate it with Meyer's quantum coordinate algebra definition used in \textsection \ref{section:quantumGrassPrelim}. This allows us to use the algebra generating set for $\OO_q(\mathrm{Gr}_{n,m})$ given in \cite{DijkStok,HK}.

The Hopf algebra $\OO_q(SU_n)$ can alternatively be presented as a Hopf $*$-subalgebra of $U_q(\frak{sl}_n)^{\circ}$, the Hopf dual of the Drinfeld--Jimbo quantum universal enveloping algebra of $\frak{sl}_n$. Explicitly, it can be defined as the Hopf subalgebra generated by the matrix coefficients of the type $1$ representations of $U_q(\frak{sl}_n)$ \cite[\textsection 2.9]{VoigtYuncken}. By construction we get a quantum generalisation of the classical Peter--Weyl decomposition
$$
\OO_q(SU_n) = \bigoplus_{\lambda \in \mathcal{P}^+} C(V_{\lambda}),
$$
where $\mathcal{P}^+$ denotes the the set of dominant integral
weights of $\frak{sl}_n$, and $C(V_{\lambda})$ denotes the coordinate algebra of $V_{\lambda}$, the irreducible representation associated to $\lambda \in \mathcal{P}^+$. 

An explicit Hopf $*$-algebra isomorphism between the two constructions can be given in the form of a non-degenerate dual pairing 
$$
\langle \cdot, \cdot \rangle: \OO_q(SU_n) \times U_q(\frak{sl}_n) \to \mathbb{C},
$$
which we consider as a Hopf algebra embedding of $\OO_q(SU_n)$ into $U_q(\frak{sl}_n)^{\circ}$. The pairing is uniquely determined by   
\begin{align} \label{eqn:dualpairing}
\langle K_i,u^j_j\rangle = q^{ \delta_{i,j-1} - \delta_{ij}}, & & 
\langle E_i,u^{i+1}_i\rangle = 1, & & \langle F_i,u^i_{i+1}\rangle = 1,
\end{align}
with all other pairings of generators being zero, see \cite[\textsection 9.4]{KSLeabh} for details. The isomorphism necessarily establishes a bijection between the cosemisimple subcoalgebras of both presentations. Explicitly it associates  $C(\lambda)$ with $C(V_{\lambda})$, where we have identified $\mathcal{P}^+$ with $\mathcal{P}_{n-1}$ in the obvious way.

Following \cite[\textsection 4]{DijkStok} we consider the Hopf subalgebra  
$$
U_q(\frak{l}_m) := \mathrm{C}\langle E_i, \, F_i, K_j \,|\, i \neq m, \, j = 1, \dots, n-1 \rangle \sseq U_q(\frak{sl}_n),
$$
where $E_i,F_i$, and $K_i$, denote the standard generators and relations of $U_q(\frak{sl}_n)$ \cite[\textsection 6]{KSLeabh}. Now the Hopf algebra $\OO_q(L_m)$ admits a dual pairing with $U_q(\frak{l}_S)$ uniquely determined by 
\begin{align*}
\langle g, \iota(X) \rangle = \langle \pi(g), X \rangle, & & \textrm{ for } X \in U_q(\frak{l}_S), \, g \in \OO_q(SU_n).
\end{align*}
Thus we see that the space of coinvariants of $\OO_q(L_m)$ coincides with the space of invariants of $U_q(\frak{l}_m)$, which is to say
\begin{align}
\OO_q(\mathrm{Gr}_{n,m}) = {}^{U_q(\frak{l}_S)}\OO_q(SU_n).
\end{align}
Hence Meyer's definition and that given in \cite[\textsection 4]{DijkStok} coincide.

The general result of \cite[Theorem 3.2]{DijkStok}, or \cite[Proposition 3.2]{HK}, restricted to the special case of the quantum Grassmannians, tells us that a set of generators of $\OO_q(\mathrm{Gr}_{n,m})$ is given by all products of the form $y\overline{y}$, where $y$ is a left highest weight element of $C(V_{\varpi_{s}})$ and $\overline{y}$ is a left lowest weight element of $C(V_{-w_0(\varpi_{s})})$, with $w_0$ denoting the longest element of the Weyl group of $\frak{sl}_n$. As is easily verified, $z_I$ is a highest weight vector of $C(V_{\varpi_{m}})$, for all $|I| = m$, and $\overline{z}_J$ is a lowest weight vector of $C(V_{-w_0(\varpi_{m})})$, for all $|J| = n-m$. Thus we see that $\OO_q(\mathrm{Gr}_{n,m})$ is generated by the set of all elements of the form $z^{IJ}$. It is now instructive to consider the Hopf subalgebra  
$$
U_q(\frak{l}^{\,\mathrm{s}}_m) := \mathrm{C}\langle E_i, \, K_i, \, F_i \,|\, i \neq m \rangle \sseq U_q(\frak{sl}_n).
$$
Using an exactly analogous argument, we can identity $\OO_q(S^{n,m})$ with the space of invariants ${}^{U_q(\frak{l}^{\mathrm{s}}_m)}\OO_q(SU_n)$.

Finally, let us consider the decomposition of the line bundles over $\OO_q(\mathrm{Gr}_{n,m})$  into irreducible left $\OO_q(SU_n)$-subcomodules, or equivalently, into irreducible right $U_q(\frak{sl}_n)$-submodules:
\begin{align*}
\EE_{\pm k} \simeq \bigoplus_{\mu \in Y} V_{\pm \mu}(k), & & \textrm{ for } k \in \mathbb{Z}_{\geq 0},
\end{align*}
where $Y$ is some subset of $\mathcal{P}^+$. As shown in \cite{DOW}, a basic weight argument confirms that each $V_{\mu}(k)$ has a highest weight element of the form $b\overline{z}^k$, for some $b \in \OO_q(\mathrm{Gr}_{n,m})$, and each $V_{\mu}(-k)$ has a highest weight element of the form $b'\!z^k$, for some $b' \in \OO_q(\mathrm{Gr}_{n,m})$. 
This fact is used in the proof of our noncommutative  Borel--Weil theorem in \textsection \ref{section:BorelWeil}.

\section{Some Additional Results}

In this appendix  we use the presentation of the Heckenberger--Kolb calculi given in \textsection \ref{section:HKRestriction} to prove some technical, yet potentially useful, results. While the results are not needed in the main body of the paper, they demonstrate the powerful tools the coquasitriangular approach gives us. We begin by producing a quantum Grassmannian representative for the coset of each basis element of $V^{(0,1)}$ and $V^{(1,0)}$. For an alternative proof, at the level of the irreducible quantum flag manifolds, see \cite[Proposition 3.6]{HKdR}.

\begin{prop}\label{nonzerolambda}
For all $(i,j) \in M^\perp \by M$, it holds that 
\begin{align} \label{eqn:zblambda}
[z^{M_{ji}M^\perp}] = \lambda_{ij} b_{ij}, & & [z^{MM^\perp_{ij}}] = \lambda_{ji} b_{ji}, 
\end{align}
for some non-zero constants $\lambda_{ij} \in \mathbb{C}$.
\end{prop}
\begin{proof}
By definition  $\{ Q_{ji} \, |\, i,j = 1, \dots, n\}$ is a dual basis to $\{b_{ij}  \, |\, i,j = 1, \dots, n\}$. Thus by Lemma \eqref{lem:QZIJ} there must exist (possibly zero) constants satisfying \ref{eqn:zblambda}. Thus the lemma would follow if we could show that these constants were non-zero. Drawing again on Lemma \ref{lem:QZIJ}, we see that  $[(z^{IJ})^+] = 0,$ for any generator $z^{IJ}$ not of the form $z^{M_{ij}M^\perp}, z^{MM^\perp_{ji}}$, for $i,j \in M \by M^\perp$. Hence,  since $\Phi\big(\Omega^1_q(\mathrm{Gr}_{n,m})\big) \in {}^H\mathrm{Mod}_0$, the elements in (\ref{eqn:zblambda}) span $V\hol \oplus V\ahol$, and so, by comparison of dimension they must all be non-zero. 
\end{proof} 

Next we give a concrete set of generators for the ideal corresponding to the Heckenberger--Kolb calculus. The argument is easily generalised to the setting of the irreducible quantum flag manifolds using \cite[\textsection 3.3.1]{HKdR}.

\begin{prop}\label{gengrass} 
The ideal corresponding to Heckenberger--Kolb calculus is generated as a right $\OO_q(\mathrm{Gr}_{n,m})$-module by the elements
\begin{enumerate}
    \item $(z^{IJ})^+  \;\;\;\qquad \mathrm{with} \; (I,J)\in T$
 \item $z^{IJ}(z^{KL})^+ \quad \mathrm{with} \; (I,J)\notin T \; \mathrm{and} \; z^{KL} \in T$,
\end{enumerate}
where we have denoted
$$
T := \left\{(M_{ij}, M^\perp),(M, M^\perp_{i'j'}),(M, M^\perp) \,|\, (i,j) \in M\times M^\perp, (i',j') \in M^\perp\times M \right\}\!.
$$
\end{prop}
\begin{proof}
As consequence of Lemma \eqref{lem:QZIJ}, the quotient of $\OO_q(\mathrm{Gr}_{n,m})^+$ by the ideal generated by the elements listed above contains the ideal corresponding to Heckenberger-Kolb calculus.
At the same time $\OO_q(\mathrm{Gr}_{n,m})^+$ is generated as right $\OO_q(\mathrm{Gr}_{n,m})$-module, by the  $(z^{IJ})^+,$ hence it follows that the dimension of the quotient of $\OO_q(\mathrm{Gr}_{n,m})^+$ by the ideal generated by the listed elements is less or equal then the dimension of the ideal corresponding to the Heckenberger-Kolb calculus, which proves the corollary.
\end{proof}

Recall next that in Lemma \ref{lem:vpmsubspaces} we showed that $V^{(0,1)}$ is an $\OO_q(SU_n)$-submodule of $\Lambda^1$. In the following proposition we produce explicit formulae for the actions of the generators  $\OO_q(SU_n)$, as well as formulae for the actions of the antipodes of the generators. This technical result will serve as a crucial ingredient in the proof that the anti-holomorphic algebra $V^{(0,\bullet)}$ is a Nichols algebra, which will appear in a forthcoming update of \cite{ROBnich}. Indeed, to match the conventions used in this forthcoming update, the result is presented in terms of the basis of $V^{(0,1)}$
$$
\{[u^i_j] \, |\, (i,j) \in M^{\perp} \times M\}.
$$
It is easily checked that this is indeed a basis of $V^{(0,1)}$ and that $[u^i_j]$ is linearly proportional to $b_{ji}$. We omit the proof which is a direct technical calculation. 
\begin{prop}
  The $\OO_q(SU_{n})$-action on~$V^{(0,1)}$ is as follows.
  For $i\neq j$ and $(r,s)\in (M\times M) \cup (M^\perp\times M^\perp)$ we have
  \begin{subequations}\label{eq:V01YDact}
    \begin{align*}
   &[u^i_j]\tl u^r_r  =  q^{-2/n}q^{\delta_{r,i}+\delta_{r,j}}[u^i_j],\\
   &[u^i_j]\tl u^r_s   =  q^{-2/n}\left(\theta(r-s)\delta_{s,i}(q-q^{-1})[u^r_j]
   + \theta(s-r)\delta_{r,j}(q-q^{-1})[u^i_s]\right)\\
   &[u^i_j]\tl S(u^r_r)  = q^{2/n}q^{-\delta_{r,i}-\delta_{r,j}}[u^i_j],\\
   &[u^i_j]\tl S(u^r_s)  =
      q^{2/n}\left( \delta_{s,i}\theta(r-s)(q^{-1}-q)[u^r_j] + \delta_{r,j}\theta(s-r)q^{2(r-s)}(q^{-1}-q)[u^i_s]
       \right)
 \end{align*}
 \end{subequations}
\end{prop}

We finish by noting that this action is a purely noncommutative phenomenon, indeed since classical forms and functions commute, the classical action of $\OO(SU_n)$ on $V^{(0,1)}$ is trivial.

\bibliographystyle{siam}

\end{document}